\documentclass[12pt,reqno]{amsart}
\usepackage{a4wide}

\usepackage{amssymb, amsmath,amsfonts, hyperref}

\usepackage{bbm}

\usepackage{graphics,graphicx,fancyhdr,color}

\newtheorem{df}{Definition}[section]
\newtheorem{corollary}[df]{Corollary}
\newtheorem{lemma}[df]{Lemma}
\newtheorem{prop}[df]{Proposition}
\newtheorem{thm}[df]{Theorem}

\newtheorem{condition}[df]{Condition}

\makeatletter \@addtoreset{equation}{section}

\newcommand{\cal}{\mathcal}
\newcommand{\bes}{\begin{displaymath}}
\newcommand{\ees}{\end{displaymath}}
\newcommand{\be}{\begin{equation}}
\newcommand{\ee}{\end{equation}}
\newcommand{\ba}{\begin{eqnarray}}
\newcommand{\ea}{\end{eqnarray}}
\newcommand{\bas}{\begin{eqnarray*}}
\newcommand{\eas}{\end{eqnarray*}}
\newcommand{\@Bbb}[1]{\ensuremath{\mathbb #1}}

\newcommand{\B}{{\@Bbb B}}
\newcommand{\C}{{\@Bbb C}}

\newcommand{\F}{{\@Bbb F}}
\renewcommand{\P}{{\mathbb P}}
\newcommand{\bbP}{{\P}}
\newcommand{\bbE}{{\mathbb E}}
\newcommand{\Q}{{\@Bbb Q}}
\newcommand{\bQ}{{\@Bbb Q}}
\newcommand{\N}{{\@Bbb N}}

\newcommand{\bbR}{{\@Bbb R}}
\newcommand{\W}{{\@Bbb W}}

\newcommand{\bbZ}{{\@Bbb Z}}
\newcommand{\bbT}{{\@Bbb T}}

\newcommand{\la}{\lambda}
\newcommand{\al}{\alpha}
\newcommand{\bt}{\beta}
\newcommand{\si}{\sigma}

\newcommand{\Om}{\Omega}
\newcommand{\om}{\omega}

\newcommand{\ep}{\varepsilon}

\newcommand{\@s}[1]{\ensuremath{\mathcal #1}}
\newcommand{\cA}{\@s A}
\newcommand{\cB}{\@s B}
\newcommand{\cC}{\@s C}
\newcommand{\cD}{\@s D}
\newcommand{\cE}{\@s E}
\newcommand{\cF}{\@s F}
\newcommand{\cG}{\@s G}
\newcommand{\cH}{\@s H}
\newcommand{\cI}{\@s I}
\newcommand{\cJ}{\@s J}

\newcommand{\cK}{\@s K}
\newcommand{\cL}{\@s L}
\newcommand{\cN}{\@s N}
\newcommand{\cM}{\@s M}
\newcommand{\cO}{\@s O}
\newcommand{\cP}{\@s P}
\newcommand{\cR}{\@s R}
\newcommand{\cS}{\@s S}
\newcommand{\cT}{\@s T}
\newcommand{\cV}{\@s V}
\newcommand{\cW}{\@s W}
\newcommand{\cX}{\@s X}
\newcommand{\cY}{\@s Y}
\newcommand{\cZ}{\@s Z}

\newcommand{\@bm}[1]{\ensuremath{\mathbf #1}}
\newcommand{\bma}{\@bm a}
\newcommand{\bmb}{\@bm b}
\newcommand{\bmc}{\@bm c}
\newcommand{\bmd}{\@bm d}

\newcommand{\bme}{\@bm e}
\newcommand{\bmf}{\@bm f}
\newcommand{\bmg}{\@bm g}
\newcommand{\bmh}{\@bm h}
\newcommand{\bmi}{\@bm i}
\newcommand{\bmj}{\@bm j}
\newcommand{\bmk}{\@bm k}
\newcommand{\bml}{\@bm l}
\newcommand{\bmm}{\@bm m}
\newcommand{\bmn}{\@bm n}
\newcommand{\bmo}{\@bm o}
\newcommand{\bmp}{\@bm p}
\newcommand{\bmq}{\@bm q}
\newcommand{\bmr}{\@bm r}
\newcommand{\bms}{\@bm s}
\newcommand{\bmt}{\@bm t}
\newcommand{\bmu}{\@bm u}
\newcommand{\bmw}{\@bm w}
\newcommand{\bmv}{\@bm v}
\newcommand{\bmx}{\@bm x}
\newcommand{\bx}{\@bm x}

\newcommand{\bmy}{\@bm y}
\newcommand{\bz}{\@bm z}

\newcommand{\by}{\@bm y}

\newcommand{\bmzero}{\@bm 0}

\newcommand{\ga}{\gamma}

\newcommand{\@g}[1]{\ensuremath{\mathfrak #1}}
\newcommand{\gA}{\@g A}
\newcommand{\gD}{\@g D}
\newcommand{\gJ}{\@g J}
\newcommand{\gF}{\@g F}
\newcommand{\gM}{\@g M}
\newcommand{\gR}{\@g R}

\newcommand{\commentout}[1]{{}}

\makeatother

\begin{document}

\title[Limit theorems for CTRW]{Limit theorems for some continuous time random walks}

\author{M. Jara}
\address{\!\!\!\!\!\!\!Ceremade, UMR CNRS 7534 \newline
  Universit\'e de Paris Dauphine,\newline
  Place du Mar\'echal De Lattre De Tassigny\newline
  75775 Paris Cedex 16 - France.
 \newline
\rm {\texttt{jara@ceremade.dauphine.fr}}}

\author{T.Komorowski}
\address{\!\!\!\!\!\!\!Institute of Mathematics, UMCS
\newline pl. Marii Curie-Sk\l odowskiej 1\newline
Lublin 20-031, Poland\newline
\!\!\!\!\!\!\!Institute of Mathematics\\
 Polish Academy of Sciences,
\newline \'{S}niadeckich 8\newline
Warsaw 00-956, Poland
\newline\rm {\texttt{
komorow@hektor.umcs.lublin.pl}}\newline
 \texttt{http://hektor.umcs.lublin.pl/\~\;\!\!komorow}}


\thanks{Work of T. K. has been partially supported by Polish
  MNiSW  grant NN 201419139 and  by EC FP6 Marie Curie ToK programme SPADE2,
  MTKD-CT-2004-014508 and Polish MNiSW SPB-M}

\begin{abstract}
In this paper we consider the scaled limit of a continuous time random walk (CTRW) based on a Markov chain
$\{X_n,\,n\ge0\}$ and two observables $\tau(\cdot)$ and $V(\cdot)$ corresponding to the renewal times and jump sizes. Assuming that these observables belong to the domains of  attraction of some stable laws we give sufficient conditions on the chain that guarantee the existence of the scaled limits for CTRW. An application of the results to a
process that arises in quantum transport theory is provided.
The results obtained in this paper generalize earlier results contained in \cite{meerschaert-1,meerschaert} and recent results of \cite{henry-straka,jurlewicz}, where $\{X_n,\,n\ge0\}$ has been a sequence of i.i.d. random variables.
\end{abstract}

\date{\today, version 3}
\maketitle

\section{Introduction}

Continuous time random walk (CTRW) has been introduced in \cite{montrol-weiss} and finds its applications in modeling of  various phenomena,  e.g. in anomalous transport (see e.g. \cite{clark-etal, gorenflo1, shlessinger1, shlessinger2, zaslavsky}) 
mathematical finance  (\cite{gorenflo,mainardi}),  or in hydrology (\cite{benson,benson1})   to mention just few examples. It can be described as a random walk subordinated to a renewal process. More precisely, suppose that $(E,$d$)$ is a Polish space with ${\cal E}$ the $\si$-algebra of its Borel subsets and $\tau:E\to(0,+\infty)$, $V:E\to \mathbb R$ are two measurable functions and $\{X_n,\,n\ge0\}$
is a Markov chain with an initial distribution $\pi$ that is stationary. Suppose also that $t_0:=0$, $t_N:=\sum_{k=0}^{N-1}\tau(X_k)$, $N\ge1$ are the renewal times. Particle jumps are given by $V(X_k)$, $k=0,1,\ldots$. Let $S_0:=0$ and $S_N:=\sum_{k=0}^{N-1}V(X_k)$.
For $t\ge0$ let $n(t):=\max[N\ge0:t_N\le t]$. We define a stochastic process describing the trajectory  of the particle performing a {\em continuous time random walk}  by
$W(t):=S_{n(t)}$, $t\ge0$. We are concerned in describing the limiting behavior of  scaled processes $\{N^{-\ga}W(Nt),\,t\ge0\}$, for an appropriate $\ga>0$, as $N\to+\infty$. In case when $\{X_n,\,n\ge0\}$ is a sequence of i.i.d. random variables this problem has been investigated in \cite{meerschaert-1}. In \cite{meerschaert} the result is  generalized to the case of triangular arrays with rowise independent random variables.. 
From Theorem 3.1 of \cite{meerschaert-1} and Theorem 2.1 of  \cite{meerschaert}   it follows  in particular that   if $\{(N^{-1/\beta}S_{[Nt]},N^{-1/\al}t_{[Nt]}),\,t\ge0\}$ converge in  law over $D([0,+\infty);\mathbb R^2)$, with the $J_1$-topology, to a L\'evy process $\{(S_t,T_t),\,t\ge0\}$, whose components have no common jumps then   $\{N^{-\al/\beta}W(Nt),\,t\ge0\}$ converges in law over $D[0,+\infty)$ with the $M_1$ topology to
\begin{equation}
\label{zeta}
\zeta_s:=S_{T^{-1}_s},\,s\ge0,
\end{equation}
where $T^{-1}_s:=\inf[t:T_t>s]$ (the first passage time) is the right inverse of the $\al$-stable subordinator $\{T_t,\,t\ge0\}$.

When common jumps of the components of  $\{(S_t,T_t),\,t\ge0\}$ are admitted with positive probability the situation is more delicate and only some partial results concerning convergence are available. In Theorem 3.4 of \cite{meerschaert-1} it has been shown that one dimensional statistics of $N^{-\al/\beta}W(Nt)$ weakly converge to the law of $\zeta_t$.

In this paper we  formulate sufficient conditions for a Markov chain $\{X_n,\,n\ge0\}$, see Theorems \ref{dependent} and \ref{dependent-no}  below, that guarantee the convergence in law of $\{N^{-\al/\beta}W(Nt),\,t\ge0\}$. 
As for the hypotheses made about the Markov chain  we assume that measure $\pi$ 
satisfies {\em the spectral gap estimate},
see Condition \ref{sg}. Moreover transition probabilities satisfy some additional regularity assumptions
see Conditions \ref{decomp} and \ref{decomp-a}.
It has been shown in \cite{jakola} that under such conditions both $\{N^{-1/\beta}S_{[Nt]},\,t\ge0\}$ and $\{N^{-1/\al}t_{[Nt]},\,t\ge0\}$ converge in  law over $D([0,+\infty)$, with the $J_1$ topology, to respective L\'evy processes. We strengthen this result and obtain the joint convergence in law of two dimensional processes $\{(N^{-1/\beta}S_{[Nt]},N^{-1/\al}t_{[Nt]}),\,t\ge0\}$ to a respective L\'evy process, see Theorems \ref{joint-conv-2}--\ref{joint-conv}. We give a sufficient condition, formulated in terms of the joint law of $(V(x),\tau(x))$ under $\pi$ (see \eqref{sing-3}), which precludes the possibility of  jumps of the
limiting L\'evy process occurring simultaneously with positive probability. The above  plus an argument from \cite{meerschaert-1} yield the convergence of $\{N^{-\al/\beta}W(Nt),\,t\ge0\}$. 
Furthermore, we prove that when the joint law of $(V(x),\tau(x))$ under $\pi$ is such that  the jumps of the
limiting L\'evy process have to occur simultaneously a.s., see condition \eqref{sing-2}, then the convergence of CTRW still holds. This is achieved by a careful analysis of the convergence of the right inverses of  $\{N^{-1/\al}t_{[Nt]},\,t\ge0\}$. We prove that the convergence of these processes holds in  a sense that allows us to control the size of the respective plateaus, see Lemma \ref{lm3.5}. This, in turn, suffices to prove the convergence of the relevant CTRW.   The limiting process is $\{\zeta_s,\,s\ge0\}$, when no common jumps are allowed for the limit  $\{(S_t,T_t),\,t\ge0\}$, see Theorem \ref{dependent}, or $\{\zeta_s^-:=S_{T_s^{-1}-},\,s\ge0\}$, see the definition of the process given in \eqref{022809} and Theorem \ref{dependent-no} below. This result is in agreement with the results obtained in the i.i.d. case in \cite{henry-straka} and \cite{jurlewicz}.

Having in mind possible applications  we  formulate the result for a counterpart of CTRW that arises when   $\{N^{-\al/\beta}S(Nt),\,t\ge0\}$ is replaced by a process  obtained by linear interpolation of its nodal points.  Finally, we apply our results to describe the limiting behavior of a jump process
$\{K_s,\,s\ge0\}$ on a one dimensional torus that arises  in  quantum transport theory,
see \eqref{lind} below. This process is the projection onto a
$0$-fiber of the solution of a  translation invariant {\em Lindblad
equation}. It possesses a unique $\si$-finite invariant measure,
absolutely continuous with respect to  Lebesgue measure, see
Proposition \ref{prop-harris} below. The dynamics of the process is  completely mixing and
its one dimensional statistics converges to a mixture of delta type measures supported on the set  $[\tau=+\infty]$, see Theorem \ref{thm1}. 
As an application of Theorem \ref{dependent} we
conclude also, see Corollary \ref{cor022101},  the convergence in law of
additive functionals of the type
$N^{-\alpha/\beta}\int_0^{Nt}V_0(K_s)ds$. In the particular case
considered  in \cite{clark-etal}  the torus is the interval
$[-\pi,\pi]$ whose endpoints are identified, $\tau(-k)=\tau(k)$ and $\tau(k)\sim |k\pm\pi/2|^{-2}$, as $|k\pm\pi/2|\ll
1$, so $\alpha=1/2$. We assume that $V_0$ is odd, i.e. $V_0(-k)=-V_0(k)$. In
addition we suppose that either $V_0(k)\sim
|k\pm\pi/2|$, as $|k\pm\pi/2|\ll
1$, and $V_0$ is bounded otherwise, or  $V_0(k)\sim
|k\pm\pi/2|^{\gamma}$, as $|k\pm\pi/2|\ll
1$ for some $\gamma>1$, and
$V_0(k)\sim  |k\pm k_0|^{-1}$ for some $k_0\not\in\{-\pi/2,\pi/2\}$.  The law
of $V(k,\tau):=V_0(k)\tau$
  belongs then to the normal domain of attraction of the Cauchy law,
  so $\beta=1/2$. We conclude therefore that the scaling
  properties of the limiting process are the same as those of the
  Brownian motion. We call  such a  process 
  {\em a fake diffusion}.

 {\bf Acknowledgement.} The authors would like to express their thanks to an anonymous  referee of the previous version of this manuscript for pointing out that  the result
 concerning the case when  the limiting process admits common jumps has been erroneously formulated.

\section{Preliminaries and statements of the main results}

\label{sec2aaa}

\subsection{A Markov chain}
\label{markov-chain}

Let $(E,$d$)$ be a  Polish metric space  and let ${\cal E}$ be its Borel
$\si$-algebra.
Assume that $\{X_n,\,n\ge0\}$ is a Markov chain with state space
$E$ and $\pi$ - the law of $X_0$ - is  \emph{invariant} and {\em
  ergodic}  for the chain. 
 We suppose that the following hypotheses are satisfies:
\commentout{Observe that an additional assumption $p(x,y)>0$ for all $x,y\in
  E$ implies that  
 $\pi$ is a unique invariant probability measure. 
}
\begin{condition}\label{sg}
{Spectral gap condition}:
  \begin{equation}
    \label{eq:8}
    \sup[\|Pf\|_{L^2(\pi)}:\,f\perp 1,\,\|f\|_{L^2(\pi)}=1]=a<1.
  \end{equation}
\end{condition}
Since $P$ is also a contraction in $L^1(\pi)$ and $L^\infty(\pi)$ we
conclude, via Riesz-Thorin interpolation theorem, that for any
$p\in[1,+\infty)$: 
\begin{equation}
  \label{eq:9}
 \|Pf\|_{L^p(\pi)}\le a^{1-|2/p-1|}\|f\|_{L^p(\pi)},
\end{equation}
for all $f\in
 L^p(\pi)$, such that $\int fd\pi=0$.

We suppose also that  the absolute continuous part of the
transition probability
function has some regularity property. Namely, we assume that:
\begin{condition}
\label{decomp}
There exist a measurable family of Borel measures $Q(x, dy)$ and a measurable, non-negative function $p(x,y)$ such that 
  \begin{equation}
P(x,dy)=P_a(x,dy)+Q(x,dy),\quad\mbox{for all }x\in E,
\label{transition}
\end{equation}
where $P_a(x,dy):=p(x,y)\pi(dy)$ and
\begin{equation}
    \label{eq:4}
    C(2):=\sup_{y\in E}\int p^2(x,y)\pi(dx)<+\infty
  \end{equation}
\end{condition}
A simple consequence of \eqref{transition} and the fact that $\pi$ is invariant is that
\begin{eqnarray}
\label{052202}
&&\int p(x,y)\pi(dy)\le 1\quad\mbox{and }\int p(y,x)\pi(dy)\le
1,\quad\forall\,x\in E.
\end{eqnarray}
Another consequence of condition \eqref{eq:4}
is that $P$ extends to a bounded operator from $L^1(\pi)$ to $L^2(\pi)$.

\subsection{The renewal process}

\label{sec3}

Suppose that $\tau:E\to [0,+\infty)$ is measurable   over
$(E,{\cal E})$ and satisfies:
\begin{condition}\label{tails}
   There exist $\alpha \in (0,1)$ and  $ c_\alpha>0$ such that
  \begin{equation}
    \label{eq:tails}
    \begin{split}
      &\lim_{\lambda\to +\infty} \lambda^\alpha\pi(\tau\ge\lambda)=
      c_\alpha
    \end{split}
  \end{equation}
\end{condition}
and there exists $t_*>0$ such that
\begin{equation*}
\label{below}
\tau(x)\ge t_*>0,\quad \forall\,x\in E.
\end{equation*}
This condition is assumed in order to avoid the issue of explosions or
accumulation points.
Furthermore, we suppose that
 the tails of $\tau$ under
the singular part are controlled by those corresponding to the absolutely continuous part uniformly with respect to the initial state, i.e. 
\begin{condition}
\label{decomp-a}
   \begin{equation}
 \label{new}
  \sup_{\la\ge0,x\in E} \frac{Q(x, [\tau\ge\lambda]) }{P_a(x, [\tau\ge\lambda]) } <+\infty. 
\end{equation}
\end{condition}
 Let $\{X_n,\,n\ge0\}$
be a Markov chain as in the previous section, $t_0:=0$ and 
\begin{equation}
\label{041901}
t_N:=\sum_{k=0}^{N-1}\tau(X_k),\quad\mbox{for $N\ge1$.}
\end{equation}
For a given $t>0$ define $n(t)$ as the unique (random) integer that satisfies the following condition
\begin{equation}
\label{011502}
  t \in[ t_{n(t)},t_{n(t)+1}).
  \end{equation}

\subsection{An observable and the CRTW process}

\label{sec2.2}

Suppose now that $V:E\to\bbR$ is measurable. 
Let $S_0:=0$
\begin{equation}
  \label{eq:101bc}
S_N:=\sum_{k=0}^{N-1}V(X_k),\quad\mbox{for $N\ge1$.}
\end{equation}
We shall assume assume that either
\begin{equation}
\label{L-2}
V\in L^2(\pi)
\quad\mbox{
and} 
\quad
\int V d\pi=0,
\end{equation}
or in case $V$ does not belong to $L^2(\pi)$ we assume that
there exist $\beta \in (0,2)$ and two 
 nonnegative constants $ c_\beta^+, c_\beta^-$ satisfying 
$ c_\beta^++
  c_\beta^->0
$ 
and 
  \begin{equation}
    \label{eq:tails:1}
    \begin{split}
      & \pi(V\ge\lambda)    =  \frac{c_\beta^+}{\lambda^{\beta}}\left(1+o(1)\right),
      \\
      &\pi(V\le-\lambda)    =  \frac{c_\beta^-}{\lambda^{\beta}}\left(1+o(1)\right),\quad \mbox{as $\lambda\to+\infty$. }
    \end{split}
  \end{equation}
 Furthermore  $V$ is supposed to be centered when $\beta\in(1,2)$. 
In analogy to condition \eqref{new} we assume that the tails of $V$ 
under the singular part of $P(x,\cdot)$ are controlled by those of the
absolutely continuous part, i.e.
   \begin{equation}
 \label{new-a}
  \sup_{\la\ge0,x\in E} \frac{Q(x, [|V|\ge\lambda]) }{ P_a(x, [|V|\ge\lambda])} <+\infty.
\end{equation}
We define the continuous time random walk (CTRW) process $W(t):=S_{n(t)}$, $t\ge0$.
Its trajectories belong to  the space of c\`adl\`ag functions $D[0,+\infty)$. To abbreviate we shall
denote this space by ${\cal D}$ in what follows.
Define also  the piecewise linear counterpart of CTRW  by
$$
\widehat W(t):=S_{n(t)}+\frac{t-t_{n(t)}}{t_{n(t)+1}-t_{n(t)}}V(X_{n(t)})\quad\mbox{for }t\in[t_{n(t)},t_{n(t)+1}).
$$ 
 In our subsequent notation we write ${\cal C}:=C[0,+\infty)$.

\subsection{Convergence to a L\'evy process}

The results presented in this section extend those of \cite{jakola} to the case of two dimensional Markov chains. They can be proved using  quite similar arguments. For the convenience of a reader we present the main points of the respective  proofs in Appendices \ref{appA} and \ref{appB}.

Suppose that  the hypotheses made in Sections \ref{markov-chain} -- \ref{sec3} hold and $K_N$ is an increasing sequence converging to infinity. 
Our immediate concern is the question of the convergence of joint processes $\{(S^{(N)}_t,T^{(N)}_t),\,t\ge0\}$, as $N\to+\infty$. Here
$T^{(N)}_t:={K_N}^{-1/\al}t_{[K_Nt]}$
 and $S^{(N)}_t:=K_N^{-1/\beta}S_{[K_Nt]}$,
when $\beta\not=1$ and 
\begin{equation}
\label{010203}
S^{(N)}_t:=K_N^{-1}S_{[K_Nt]}-v_Nt,
\end{equation}
when $\beta=1$ and $v_N:=\int V1_{[|V|\le K_N]}d\pi$.
Let
 \begin{equation}
\label{joint-1}
\psi_\beta(\xi):=\int_{\bbR}e_\beta(\xi,\lambda)\nu_{\beta}(d\la),
\end{equation}
where 
$$
e_\beta(\xi,\lambda):=
\left\{
\begin{array}{ll}
e^{i\lambda
  \xi}-1,&\beta\in(0,1),\\
e^{i\xi \la}-1-i\xi\la1_{[-1,1]}(\la),&\beta=1,\\
e^{i\lambda \xi}-1-i\lambda\xi,& \beta\in(1,2)\\
\end{array}
\right.
$$
and
$\nu_{\beta}(d\la):=\beta c_\beta(\la)|\la|^{-1-\beta}d\la$. Here $c_{\beta}(\lambda)$ equals $c_\beta^+$ for $\lambda>0$ and $c_\beta^-$ for $\la<0$. 
Consider a L\'evy process $\{(S_t,T_t),\,t\ge0\}$ given by
$\bbE e^{i\xi_1 S_t+i\xi_2 T_t }=e^{t\psi(\xi_1,\xi_2)}$,
where 
\begin{eqnarray}
\label{joint-1bbc}
&&\psi(\xi_1,\xi_2):=\psi_\beta(\xi_1)+\psi_\al(\xi_2)\\
&&
=\int e_{\al,\bt}(\xi_1,\xi_2,\lambda_1,\lambda_2)\nu_*(d\la_1,d\la_2)
,\quad (\xi_1,\xi_2)\in\mathbb R^2,\nonumber
\end{eqnarray}
where the coefficients $c_\al^+,c_\al^-$ appearing in the respective definition of $\psi_\al(\cdot)$ are equal then to $c_\al$ (see \eqref{eq:tails})  and $0$ respectively,
\begin{equation}
\label{011902}
e_{\al,\bt}(\xi_1,\xi_2,\lambda_1,\lambda_2):=
\left\{
\begin{array}{ll}
e^{i(\lambda_1
  \xi_1+\lambda_2\xi_2)}-1,&\beta\in(0,1),\\
e^{i(\lambda_1
  \xi_1+\lambda_2\xi_2)}-1-i\xi_1\la_11_{[-1,1]}(\la_1),&\beta=1,\\
e^{i(\lambda_1
  \xi_1+\lambda_2\xi_2)}-1-i\lambda_1\xi_1,& \beta\in(1,2)\\
\end{array}
\right.
\end{equation}
and 
\begin{equation}
\label{021902}
\nu_*(d\la_1,d\la_2):=\nu_{\beta}(d\la_1)\delta_0(d\la_2)+\delta_0(d\la_1)\nu_{\al}(d\la_2).
\end{equation}
In our first result we adopt a hypothesis that $\tau(x)$ and
$|V(x)|$ cannot be large together. Namely, we assume that
\begin{equation}
\label{sing-3}
\pi\left[\tau\ge \lambda,\, |V|\ge \lambda\right] \le \frac{C_*}{\lambda^{\gamma}},\quad\lambda>0
\end{equation}
for some $C_*>0$ and $\gamma>\alpha\vee \bt$. 
\begin{thm}
\label{joint-conv-2}
Suppose that the assumptions made in \eqref{eq:tails:1} and  \eqref{sing-3} hold and 
$V$ is centered
when $\beta\in(1,2)$. Then, the following hold:

i) if $\beta\not=1$ then the joint laws of
 $\{(S^{(N)}_t,T^{(N)}_t),\,t\ge0\})$ converge in law, as
 $N\to+\infty$, to $\{(S_t,T_t),\,t\ge0\})$ on
${\cal D}_2:=D([0,+\infty),\bbR^2)$ with the $J_1$-topology.

ii) if $\beta=1$, we assume that for some $\beta'>1$ 
\begin{equation}
\label{020411b}
\sup_{N\ge1}\|PV_N\|_{L^{\beta'}(\pi)}<+\infty,
\end{equation}
where $V_N:=V1_{[|V|\le K_N]}$. 
Then, the conclusion of part i) holds also in this case with the  modification of the definition of $S^{(N)}_t$ given in \eqref{010203}.
\end{thm}

In our next result we allow the jumps of the components of 
$\{(S_t^{(N)},T_t^{(N)}),\,t\ge0\}$ to occur at the
same time. More specifically,
let 
\begin{equation}
\label{011802} 
\rho(\lambda):=C_{\al,\bt}|\lambda|^{\bt/\alpha},\quad\lambda\in \mathbb R
 \end{equation} 
where $C_{\al,\bt}:=c_\al(c_{\bt}^-+c_{\bt}^+)^{-1}$.
Suppose that for some $C_*>0$ and $\gamma>\alpha$ we have
\begin{equation}
\label{sing-2}
\pi\left[|\tau-\rho\circ V|\ge \lambda\right] \le \frac{C_*}{\lambda^{\gamma}},\quad\forall\,\lambda>0.
\end{equation}
Consider now a L\'evy process $\{(S_t,T_t),\,t\ge0\}$ such that
\begin{equation}
\label{joint-1bb}
\psi(\xi_1,\xi_2):=\int_{\bbR^2}e_{\al,\bt}(\xi_1,\xi_2,\lambda_1,\lambda_2)\nu_*(d\la_1,d\la_2),
\end{equation}
where 
$e_{\al,\bt}(\xi_1,\xi_2,\lambda_1,\lambda_2)$ is given by \eqref{011902}
and
\begin{equation}
\label{031902}
\nu_*(d\lambda_1,d\lambda_2):=\delta_0(\lambda_2-\rho(\lambda_1))\nu_{\beta}(d\lambda_1)d\lambda_2.
\end{equation} 
 \begin{thm}
\label{joint-conv-1}
Suppose  that \eqref{sing-2} is in force. Then, the convergence statements analogous to the one made in parts
i), ii) of Theorem \ref{joint-conv-2} still hold. The only difference is that  the limiting   L\'evy process is   described by the exponent given in \eqref{joint-1bb}. 
\end{thm}

Finally, when $V\in L^2(\pi)$, i.e. \eqref{L-2} holds, we have  the following.
\begin{thm}
\label{joint-conv}
The  laws of  $\{(S^{(N)}_t,T^{(N)}_t),\,t\ge0\})$ converge, as
$N\to+\infty$,   over ${\cal D}\times {\cal D}$ with the product  of the uniform and $J_1$ topologies to the joint law of independent L\'evy
processes: $\{(S_t,T_t),\,t\ge0\}$. The first component is a zero mean Brownian motion and the L\'evy exponent of the second component equals $\psi_\al(\xi)$.
\end{thm}

\subsection{Convergence of continuous time random walks}

Our first result, concerning the convergence of CTRW, is contained in the following.
\begin{thm}
\label{dependent}
Under the assumptions of either Theorem \ref{joint-conv-2},  or \ref{joint-conv} the processes \linebreak
$\{N^{-\al/\beta}W(Nt),\,t\ge0\}$ converge in law
in the $M_1$ topology  of ${\cal D}$, as $N\to+\infty$, 
to 
$\{\zeta_t:=S_{s(t)},\,t\ge0\}$, where $\{(S_t,T_t),\,t\ge0\}$ is an appropriate
L\'evy process and
$\{s(t),\,t\ge0\}$ is the right inverse 
of $\{T_t,\,t\ge0\}$. The result also holds when process $W(t)$
is replaced by the linear interpolation process $\widehat W(t)$. In the latter case, under the assumptions of Theorem \ref{joint-conv} the convergence in law holds over ${\cal C}$. 
\end{thm}
{\bf Remark 1.} If $\{(S_t,T_t),\,t\ge0\}$ is such that its first component is a Brownian motion then the components of the process are independent, see Theorem \ref{joint-conv}.  In that case  $\{\zeta_s=S_{T^{-1}_s},\,s\ge0\}$ is called a {\em Mittag-Leffler process}. It is non-Markovian and arises as a limit of an appropriately scaled additive
functional of a Markov process, whose resolvent, applied at the observable, obeys the power law at  the bottom of the spectrum of the generator,
see \cite{athreya,darling-kac}. We refer the reader to
e.g. \cite{hopfner} and the references therein for
an extensive review of  the results concerning
this particular case.  

Given the L\'evy process   $\{(S_t,T_t),\,t\ge0\}$ as described  in Theorem   \ref{joint-conv-1} 
define
\begin{equation}
\label{022809}
\zeta_t^-=\lim_{N\to+\infty}S_{s(t)-1/N},\quad t\ge0.
\end{equation}
The limit is understood   almost surely  in the $J_1$ topology of ${\cal D}$. Here $S_t=0$ when $t<0$. Observe that although the notation suggests otherwise the process $\{\zeta_t^-,\,t\ge0\}$ is  c\`adl\`ag, as a limit of c\`adl\`ag processes in the $J_1$ topology..
\begin{thm}
\label{dependent-no}
Under the assumptions of  Theorem   \ref{joint-conv-1} the processes
$\{N^{-\al/\beta}W(Nt),\,t\ge0\}$ converge in law
in the $J_1$ topology  of ${\cal D}$, as $N\to+\infty$, 
to 
$\{\zeta_t^-,\,t\ge0\}$, defined above. 
\end{thm}

{\bf Remark 2.}
We point out here that the limiting processws
 described in Theorems \ref{dependent} and \ref{dependent-no} have  a scale invariance property. Namely, the laws of
$\{\zeta_{at},\,t\ge0\}$ and that of  $\{a^{\alpha/\bt}\zeta_{t},\,t\ge0\}$
are identical for each $a>0$. 
The same scaling invariance concerns also the process $\{\zeta_{t}^-,\,t\ge0\}$
This remark follows easily from the fact
that in the cases considered in both theorems the L\'evy processes
$\{(S_{t a^{\al}},T_{t a^{\al}}),\,t\ge0\}$ and
$\{(a^{\al/\bt}S_{t},aT_{t}),\,t\ge0\}$ have identical L\'evy
exponents. Thus, the joint laws of $\{(S_{t a^{\al}},T_{t a^{\al}}),\,t\ge0\}$ and those of
$\{(a^{\al/\bt}S_{t},aT_{t}),\,t\ge0\}$ over ${\cal
  D}_2$ are identical. This in turn implies  easily the scale invariance property.

\section{The proof of Theorem \ref{dependent}}

\subsection{The case when jumps cannot occur together}
Here we assume that the sets of discontinuity points for the components of the limiting process $\{S_t,\,t\ge0\}$ and  $\{T_t,\,t\ge0\}$ are a.s. disjoint, i.e. that either the assumptions of Theorem \ref{joint-conv-2}, or \ref{joint-conv} hold. Then, the weak convergence of $\{N^{-\al/\beta}W(Nt),\,t\ge0\}$ can be proven in exactly the same way
 as in Theorem 3.1 of \cite{meerschaert-1}. 
We
 only  show the convergence of the linear interpolation process $N^{-\al/\beta}\widehat W(Nt)$. 
 
 Assume first that the assumptions of Theorem \ref{joint-conv-2} hold.
 For a given $t>0$ recall that $n(Nt)$ is the (random) integer given by \eqref{011502}. Let $K_N:=N^{\al}$.
We define
$
 S^{(N)}_t:=K_N^{-1/\beta}\sum_{n=0}^{[K_Nt]-1}V(X_n),
$
$ T^{(N)}_t:=K^{-1/\al}_N t_{[K_Nt]}$
 and $ s_N(t)$   the right-continuous inverse of $\{ T^{(N)}_t,\,t\ge0\}$, i.e. $ s_N(t):=\inf[u: T^{(N)}_u>t]$, and $s_{N}^*(t):=s_N(t)-N^{-\al}=\max[u: T^{(N)}_u\le t]$. 
Denote by $\widehat S_t^{(N)}$ the process whose paths are obtained by
the linear interpolation between the points $(mK_N^{-1}, S_{mK_N^{-1}})$, where
$m\ge0$ is an integer. 
Then $N^{-\al/\beta}\widehat W(Nt)=\widehat S^{(N)}_{s_{N}^*(t)}$.
The following result allows us to
replace the first coordinate process in the statement of Theorem
\ref{joint-conv-2} by its linear interpolation. 
\begin{lemma}
\label{lm-aux-1}
Under the assumptions of Theorem \ref{joint-conv-2} the processes
$\{(\widehat S^{(N)}_t, T^{(N)}_t),\,t\ge0\}$ converge in law, as
$N\to+\infty$, over  ${\cal D}\times {\cal D}$ with the product of the $M_1$ topologies to the  L\'evy process: $\{(S_t,T_t),\,t\ge0\}$  as in the statement of Theorem \ref{joint-conv-2}.
\end{lemma}
Before demonstrating the lemma we show how to use it to finish the proof of the theorem.  By Skorochod's embedding theorem we define a family of processes
 $\{(U^{(N)}_t,V^{(N)}_t),\,t\ge0\}$ such that
\begin{itemize}
 \item[1)] the law of  $\{(U^{(N)}_t,V^{(N)}_t),\,t\ge0\})$
is identical with that of $\{ (\widehat S^{(N)}_t,  T^{(N)}_t),\,t\ge0\})$ for each $N\ge1$,
\item[2)]  $\{(U^{(N)}_t,V^{(N)}_t),\,t\ge0\}$
  converges a.s., in ${\cal D}\times {\cal D}$ equipped with  the
 product $M_1$ topology,
to $\{(S_t, T_t),\,t\ge0\})$. Here the limiting process is as in Theorem \ref{joint-conv-2}.
\end{itemize}
Suppose that
 $\{u_N(t),\,t\ge0\}$ is the right inverse of $\{V^{(N)}_t,\,t\ge0\}$ and $u_N^*(t):=u_N(t)-1/N^{\al}$.
The law of 
$
\{Y^{(N)}_t:=U^{(N)}_{u_N^*(t)},\,t\ge0\}$ coincides with that of $N^{-\al/\beta}\widehat W(Nt)$.
Moreover,  both $
\{U^{(N)}_t,\,t\ge0\}$ and  $\{u_N^*(t),\,t\ge0\}$ converge
a.s., as $N\to+\infty$, in the $M_1$ topology to $\{S_t,\,t\ge0\}$ and $\{s(t),\,t\ge0\}$  (the right  inverse of $\{T_t,\,t\ge0\}$) respectively. In fact, since $\{s(t),\,t\ge0\}$ is a.s. continuous, the latter sequence converges in the uniform topology.   Theorem 13.2.4 of \cite{whitt}
implies therefore that $\{Y^{(N)}_t,\,t\ge0\}$ converge in the $M_1$ topology to
$\{Y_t,\,t\ge0\}$ a.s., provided that the sets of dicontinuities of  $\{S_t,\,t\ge0\}$ and $\{
 T_t,\,t\ge0\}$ are a.s. disjoint. This however is a simple
 consequence of the independence of these processes.

\subsection*{The proof of Lemma \ref{lm-aux-1}}
Suppose $T>0$. Recall how  the $M_1$ topology on $D[0,T]$ can be metrized, see
\cite{whitt}, p. 476 for details. For a given  $X\in
D[0,T]$
we define by $\Gamma_X$ the graph of $X$, i.e. the subset of $\bbR^2$
given by 
$$
\Gamma_X:=[(t,z):t\in[0,T],\, z=c X(t-)+(1-c)X(t)],\mbox{ for some
}c\in[0,1]].
$$ 
On $\Gamma_X$ we define an order by letting $(t_1,z_1)\le(t_2,z_2)$
iff
$t_1< t_2$, or $t_1=t_2$ and $|X(t_1-)-z_1|\le |X(t_1-)-z_2|$. Denote
by $\Pi(X)$ the set of all continuous mappings $\gamma=(\gamma^{(1)},\gamma^{(2)}):[0,1]\to\Gamma_X$ that
are non-decreasing, i.e. $t_1\le t_2$ implies that $\gamma(t_1)\le
\gamma(t_2)$.
The metric d$(\cdot,\cdot)$ is defined as follows:
$$
 \mbox{d}(X_1,X_2):=\inf[\|\gamma^{(1)}_1-\gamma^{(1)}_2\|_\infty\vee \|\gamma^{(2)}_1-\gamma^{(2)}_2\|_\infty,\,\gamma_i=(\gamma_i^{(1)},\gamma_i^{(2)})\in \Pi(X_i),\,i=1,2].
$$
This metric provides a metrization of the $M_1$ topology, see
\cite{whitt}, Theorem 13.2.1.

For any $\gamma_1\in
\Gamma_{ S^{(N)}_\cdot}$ we  define  $\gamma_2\in
\Gamma_{\widehat S^{(N)}_\cdot}$ as follows. Suppose $\gamma_1(t)$
belongs to the graph corresponding to $(t, S^{(N)}_t)$, $t\in
[mK_N^{-1},(m+1)K_N^{-1})$ for an integer $m\ge0$. Let $\gamma_2(t)$
be the nearest neighbor projection of $\gamma_1(t)$ onto the segment
joining $(mK_N^{-1}, S^{(N)}_{ mK_N^{-1}})$ with $((m+1)K_N^{-1}, S^{(N)}_{(m+1)K_N^{-1}-})$. One
can use these two parametrizations to estimate the distance d as follows
$$
 \mbox{d}(S^{(N)}_\cdot,\widehat S^{(N)}_\cdot)\le C K_N^{-1/2(1+1/\beta)}\max_{0\le n\le [T K_N]}|V(X_n)|^{1/2}
$$
for some deterministic constant $C>0$ independent of $N$.
Hence for any $\eta>0$
we get
\begin{equation}
\label{042102}
\bbP\left[\mbox{d}( S^{(N)}_\cdot,\widehat S^{(N)}_\cdot)\ge
  \eta\right]\le CK_N^{-\bt}\to0,
\end{equation}
as $N\to+\infty$. The lemma is then a consequence of  \eqref{042102} and Theorem \ref{joint-conv-2}\qed

When, on the other hand, the assumptions of Theorem \ref{joint-conv} hold we can conclude that
for each $T,\eta>0$
\begin{equation}
\label{011602}
\lim_{N\to+\infty}\bbP\left[\sup_{t\in[0,T]}| S^{(N)}_t-\widehat S^{(N)}_t|\ge
  \eta\right]=0,
\end{equation}
which implies the weak convergence of  the linear interpolation process $\{N^{-\al/\beta}\widehat W(Nt),\,t\ge0\}$ over ${\cal C}$ to a Mittag-Leffler process $\{\zeta_t,\,t\ge0\}$.
To show \eqref{011602} note that the expression under the limit on its left hand side can be estimated from above by
\begin{eqnarray}
\label{021602}
&&\bbP\left[\max_{0\le n\le K_NT}|V(X_n)|\ge
  K_N^{1/2}\eta\right]\le K_NT\pi\left(|V|\ge
  K_N^{1/2}\eta\right)\\
  &&\le \frac{K_N T}{(\eta K_N^{1/2})^2}\int_{[|V|\ge
  K_N^{-1/2}\eta]}V^2d\pi= \frac{T}{\eta^2}\int_{[|V|\ge
  K_N^{-1/2}\eta]}V^2d\pi\to0,\nonumber
\end{eqnarray}
as $N\to+\infty$.

\subsection{The case when jumps occur together} We assume the hypotheses of Theorem \ref{joint-conv-1} and 
\label{sec4.3}
 admit $K_N:=N^{\al}$, as in the previous section. Using  Skorochod's embedding theorem define a family of
processes
$\{(U^{(N)}_t,V^{(N)}_t),\,t\ge0\}$ such that:
\begin{itemize}
 \item[1)] the law of $\{(U^{(N)}_t,V^{(N)}_t),\,t\ge0\}$
coincides with that of $\{ ( S^{(N)}_t,T^{(N)}_t),\,t\ge0\}$ for each $N\ge1$,
\item[2)]  $\{(U^{(N)}_t,V^{(N)}_t),\,t\ge0\}$
  converges a.s., in the
  $J_1$ topology of ${\cal D}_2$,
to $\{(S_t,T_t),\,t\ge0\}$. The latter process is as in the
statement of Theorem \ref{joint-conv-1}. The above means that for any $L>0$ one can
find a sequence $\{\lambda_n;n \geq 1\}$  of increasing homeomorphisms
in $[0,L]$ such that
$\lambda_N(0)=0$, $\lambda_N(L)=L$ and 
\begin{equation}
\label{020504}
 \sup_{t\in[0,L]}|\lambda_N(t)-t|\to 0,
\end{equation} 
\begin{equation}
\label{030504}
 \sup_{t\in[0,L]}|U^{(N)}_{\lambda_N(t)}-S_t|\to
 0\quad\mbox{ and}\quad \sup_{t\in[0,L]}|V^{(N)}_{\lambda_N(t)}-T_t|\to 0,
\end{equation} 
as $N\to+\infty$.
\end{itemize}
Let $u_N(t)$, $s(t)$ be the right inverses of  $V_t^{(N)}$, $T_t$ respectively. On the other hand, if $u_N^*(t):=\max[s:U^{(N)}_s\le t]$ then
$u_N^*(t)=u_N(t)-1/N^{\al}$. 
Observe that the CTRW $\{N^{\al/\beta}W(Nt),\,t\ge0\}$ has the same law as $\{U^{(N)}_{u_N^*(t)},\,t\ge0\}$.

We show the following.
\begin{thm}
\label{lm-x-1}
Under the assumptions of Theorem \ref{joint-conv-1} 
the processes $\{U^{(N)}_{u_N^*(t)},\,t\ge0\}$
converge in law over ${\cal D}$, with the $J_1$ topology, as $N\to+\infty$, 
to  $\{\zeta_t^-,\,t\ge0\}$ defined in \eqref{022809}.
\end{thm}

{\em Proof.} 
The proof relies on a careful analysis of  convergence of processes $u_N(\cdot)$ to $s(\cdot)$. We know that they converge uniformly on any compact interval. In fact as we show in Lemma \ref{lm3.5} below for $t$ in a   plateau of $s(\cdot)$ of a fixed size we have 
$\la_N^{-1}\circ u_N(t)=s(t)$ for a sufficiently large $N$. Matching plateaus of $s(\cdot)$ with those of $\la_N^{-1}\circ u_N(\cdot)$ we define homeomorphisms $\Lambda_N$ such that 
$\lim_{N\to+\infty}U^{(N)}_{u_N\circ \Lambda_N(t)}= S_{s(t)}$ uniformly on the set of plateaus of $s(\cdot)$ of a fixed size. To show that this convergence extends also to the entire $[0,T_L]$ we use the fact that the  plateaus of $s(\cdot)$ correspond also to the  jumps of $\{S_t,\,t\ge0\}$, see Lemma \ref{uniform}. This is due to the fact that  jumps of $S_t$ and $T_t$  are matched by function $\rho(\cdot)$.  Therefore, outside the large size plateaus of $s(\cdot)$ the trajectory  $\{S_{s(t)},\,t\ge0\}$ cannot suffer large jumps.  Since  $|u_N\circ \Lambda_N(t)-s(t)|$ is small for sufficiently large $N$ we can easily conclude  that $|U^{(N)}_{u_N\circ \Lambda_N}(t)- S_{s(t)}|$ is also small. This fact implies also that $\lim_{N\to+\infty}|U^{(N)}_{u_N^*\circ \Lambda_N}(t)- S_{s(t)-1/N^{\al}}|=0$ and the theorem follows.

 To start the rigorous proof  we need some auxiliary results.
For each $\delta>0$, we denote by $A_\delta = A_\delta(s(\cdot))$ the set of
 "plateau points'' of $s(t)$ of size at least $\delta$, i.e. 
\begin{equation}
\label{plateau}
\mbox{$t \in A_\delta$ iff
 $s(t)$ is constant in  the interval $(t-\delta,t+\delta)\cap
 [0,T_{L}]$.}
\end{equation}
Obviously $A_\delta\subset A_{\delta'}$ if $\delta'\le \delta$.
\begin{lemma}
\label{lm3.5}
For a fixed $L>0$ the sequence $\{u_N(\cdot),\,N\ge1\}$ converges to $s(\cdot)$ in the following sense:
\begin{itemize}
\item[i)] $\|\lambda_N^{-1} \circ u_N -s\|_\infty\to0$, a.s. as
  $N\to+\infty$, where the supremum norm is taken over $[0,L]$,
\item[ii)] there exists a decreasing sequence $d_N\to 0$, as $N\to+\infty$,  such
  that 
$\lambda_n^{-1}
  \circ u_n(t)= s(t)$ for $t \in A_{d_N}$, $n\ge N$ and
  all $N\ge1$.
\end{itemize}
\end{lemma}

{\em Proof.}
Let
\begin{equation}
\label{060505}
a_N
:=\max[\|V^{(N)} \circ \lambda_N -T\|_\infty, \|U^{(N)} \circ
  \lambda_N -S\|_\infty, \|\lambda_N
  -\text{id}\|_\infty,\|\lambda_N^{-1} -\text{id}\|_\infty]
\end{equation}
and
\begin{equation}
\label{a-n}
d_N
:=\sup_{n\ge N}[\|V^{(n)} \circ \lambda_n -T\|_\infty, \|U^{(n)} \circ
  \lambda_n -S\|_\infty, \|\lambda_n
  -\text{id}\|_\infty,\|\lambda_n^{-1} -\text{id}\|_\infty].
\end{equation}
From \eqref{020504} and \eqref{030504} we have $d_N\to0$, thus also,  $a_N\to0$, as $N\to+\infty$.
Since each $u_N(t)$  is increasing and $s(t)$ is a.s. continuous in order to prove i) it is enough to show that a.s.
$$
\lim_{N\to+\infty}u_N(t)=s(t),\quad\forall\,t\ge0.
$$ 
We prove first that the convergence holds a.s. for any fixed $t$. We claim that then
\begin{equation}
\label{s-1}
 u_N(t) \leq  s(t+a_N)+a_N,\quad \mbox{a.s.}
\end{equation}
Indeed,
let $\delta>0$ and
$ s_*:=s(t+a_N)+a_N+\delta$. Then $ s(t+a_N)=s_*-a_N-\delta<
\lambda_N^{-1}(s_*) $, by definition \eqref{060505}. From the definition of the right inverse $t+a_N \le T_{s(t+a_N)}$.
Since $\{T_s,\,s\ge0\}$ is a.s. strictly increasing (\cite{sato}, Theorem 21.3, p. 136) 
$$
V^{(N)}_{s_*}\ge T_{\lambda_N^{-1}(s_*)}-a_N> T_{s(t+a_N)}-a_N\ge t,
$$ so
 $$
 u_N(t)\le s_*=s(t+a_N)+a_N+\delta.
 $$ 
 Since $\delta>0$ was arbitrary  \eqref{s-1} follows.

Likewise, we prove that $u_N(t) \geq s(t-a_N)-a_N$, a.s. and, as a
result, we
conclude that there exists $\mathbb  D$, a dense subset of $[0,+\infty)$, such that
$\lim_{N\to+\infty}u_N(t)=s(t)$ for $t\in \mathbb D$ a.s.
\commentout{

Suppose that $\delta>0$ is arbitrary. We claim that
$$
u_N(t)+a_N> s(t-a_N)-\delta.
$$
Indeed, assume that
$
u_N(t)+a_N\le s(t-a_N)-\delta,
$
then $\lambda_N^{-1}(u_N(t))\le s(t-a_N)-\delta$ and 
$$
t-a_N\le V^{(N)}_{u_N(t)}-a_N\le T_{\lambda_N^{-1}(u_N(t))}\le T_{s(t-a_N)-\delta}<T_{s(t-a_N)-\delta/2}
$$
But this is impossible because then $s(t-a_N)\le s(t-a_N)-\delta/2$.

}
  Since all the functions $u_N$ are increasing the convergence can be easily extrapolated to the entire  $[0,+\infty)$. This ends the proof of part i)  of the lemma.

For
 any $t \in A_{d_N}$ and $n\ge N$ we have: 
\begin{eqnarray*}
&&
T_{ s(t)-} \leq t-d_N \le t+d_N\leq T_{ s(t)},\\
&&
|V^{(n)}_{\lambda_n\circ s(t)-} - T_{ s(t)-}|\leq d_N
\quad\mbox{
and} 
\quad
|V^{(n)}_{\lambda_n\circ  s(t)} -T_{ s(t)}|\leq d_N.
\end{eqnarray*}
 Therefore, $V^{(n)}_{\lambda_n\circ s(t)-} \leq t \leq
 V^{(n)}_{\lambda_n\circ  s(t)}$, which proves that $u_n(t)=\lambda_n \circ  s(t)$. 
This ends the proof of part ii) and thus completes the demonstration of the lemma.$\Box$

For each $t\ge0$  define $\Delta S_t =S_t-S_{t-}$.

\begin{lemma} 
\label{uniform}
Let $\{d_m,\,m\ge1\}$ be as in the statement of Lemma \ref{lm3.5}.
Under the assumptions of Theorem \ref{joint-conv-1} we have
\begin{equation}
\label{050505}
\lim_{m \to \infty} \sup_{t \notin s(A_{d_m})} |\Delta S_t| =0\quad \text{in probability}
\end{equation}
\end{lemma}
\proof
Let 
$$
B_m:=[\exists\,t\in[0,L]:\,|\Delta S_t|\ge
\rho^{-1}(4d_m)\,\text{ and }\Delta T_t\le d_m],
$$
where $\rho(\cdot)$ is given by \eqref{011802}.
We show that 
\begin{equation}
\label{030705}
\bbP[B_m]=0,\quad\forall\,m\ge1.
\end{equation} 
Suppose first that $\beta<1$.
Consider the  jump process $\{Z_t^{(r)}:=(S_t^{(r)},T_t^{(r)}),\,t\ge0\}$ corresponding to the jump 
measure
$$
\nu_*^{(r)}(d\lambda_1,d\lambda_2):=1_{B_r^c(0)}(\lambda_1,\lambda_2)
\nu_*(d\lambda_1,d\lambda_2). 
$$

Let $z:=\nu_*^{(r)}(\bbR^2)$.
This process can be realized as follows: $Z_t^{(r)}=
Z_{N(t)}^{(r)}$, where $Z_{n}^{(r)}$ is a sum of $n$ independent
random variables distributed according to
$z^{-1}\nu_*^{(r)}(d\lambda_1,d\lambda_2)$
and $N(t)$ is an independent of it Poisson process with intensity
$z$. Since the jumps of $\{Z_t^{(r)},\,t\ge0\}$ are vectors whose coordinates belong to the support of $\nu_*^{(r)}(d\lambda_1,d\lambda_2)$ 
(contained in the curve $\{(\lambda,\rho(\lambda)),\lambda>0\}$) we have
\begin{equation}
\label{010705}
\bbP[B_m^{(r)}]=0,
\end{equation} 
where 
$$
B_m^{(r)}:=[\exists\,t\in[0,L]:\,|\Delta S_t^{(r)}|\ge
\rho^{-1}(3d_m)\,\text{ and }\Delta T_t^{(r)}\le 2d_m].
$$

Let $\{Z_t:=(S_t,
T_t),\,t\ge0\}$. It is well known, see
e.g. \cite{breiman} Theorem 14.27, that
\begin{equation}
\label{020705}
\lim_{r\to0+}\sup_{t\in[0,L]}|Z_t^{(r)}-Z_t|=0,\quad
\text{in probability}.
\end{equation} 
Combining \eqref{010705} and \eqref{020705} we conclude
\eqref{030705}. Thus \eqref{050505} follows.

The case when $\beta\in[1,2)$ can be concluded similarly. However,
then the approximating processes should be of the form $\{Z_t^{(r)}-c^{(r)}t, \,t\ge0\}$ for some
$c^{(r)}=(c^{(r)}_1,0)$, where, in general, $c^{(r)}_1$ may diverge, as $r\to0+$. 
\qed

\subsection*{Proof of Theorem \ref{lm-x-1}}


%
{\bf Step 1.} First, we show that $
\lim_{N\to+\infty}U^{(N)}_{u_N(t)}= S_{s(t)}$ in the $J_1$ topology.
Writing $U^{(N)} _{u_N (t)} =U^{(N)}
_{\lambda_N\lambda_N^{-1}u_N (t)} $, we notice that, in light of \eqref{030504}, it
is enough to show convergence in the $J_1$-Skorohod topology of $S_{ \lambda_N^{-1}u_N (t) }$ to $S_{
  s(t)}$. For any $L>0$ we exhibit increasing homeomorphisms
$\Lambda_N: [0,T_L] \to [0,V_L^{(N)}]$, $N\ge1$ such that  
\begin{equation}
\label{020505}
\lim_{N\to+\infty}
S_{ \sigma_N(t) } =S_{
  s (t)}
\end{equation}
and 
\begin{equation}
\label{030505}
\lim_{N\to+\infty}\Lambda_N(t) =t,
\end{equation}
uniformly on $[0, T_L]$. 
Here
\begin{equation}
\sigma_N(t):=\lambda_N^{-1}\circ u_N\circ \Lambda_N (t).
\label{012809}
\end{equation}
We  can conclude from the above argument and from \eqref{030504} that
$\lim_{N\to+\infty}[U^{(N)}_{ u_N\circ\Lambda_N(t)}-S_{s(t)}]=0$ uniformly on compact intervals.

We display now the construction of the homeomorphisms $\Lambda_N(t)$ that saisfy \eqref{020505} and \eqref{030505}. Suppose that 
$\{d_n,\,n\ge1\}$ is strictly decreasing sequence as in the statement of Lemma \ref{lm3.5}. Let $\{\ell_k,\,k \ge1\}$ be an increasing sequence of positive
 integers such that ${\cal S}_{k} :=\{t_1,\ldots,t_{\ell_k}\}= s(A_{d_k})$. Then,
 $$
s^{-1}({\cal S}_{k} )=\bigcup_{i=1}^{\ell_k}[T_{t_i-},T_{t_i}]\supset A_{d_k}.
$$
Intervals $[T_{t_i-},T_{t_i}]$
 are the  plateaus of $s(t)$. They are mutually disjoint and each is of length
greater than $2d_k$. The complement of  $s^{-1}({\cal S}_{k} )$ in $[0,T_L]$ is an open set that is a union of a
finite number of open intervals (relative to $[0,T_L]$). Let
$\kappa_{k}$ be the minimum of the lengths of  these
intervals. Of course $\kappa_k$ decreases to $0$, as $k\to+\infty$.
Let $\{m_k,\,k\ge1\} $ be an increasing sequence of positive integers such that $d_N <
\min[\kappa_{k}/2,d_k-d_{k+1}]$ for all $N \geq m_k$.
 Recall that then both $|V^{(N)}_{
\lambda_{N} (t_i)-}-T_{t_i-}| $ and $|V^{(N)}_{
\lambda_{N} (t_i)}- T_{t_i}| $, $i=1,\ldots,\ell_k$ are  less than, or equal to
$d_{m_k}$ for $N\ge m_k$. 
Therefore,  for each such
$N$ the intervals $[V^{(N)}_{\lambda_{N}(t_i)-},V^{(N)}_{\lambda_{N} (t_i)}]$ (plateaus
of $u_{N}$) are mutually disjoint for different $i=1,\ldots,\ell_k$. 

We say that the interval $[c,d]$ follows
$[a,b]$ if $c> b$. Let us
take $i$, $j$ such that their corresponding plateaus $[
T_{t_i-}, T_{t_i}]$ and $[ T_{t_j-}, T_{t_j}]$ are
consecutive (in this order). Then $[
V^{(N)}_{\lambda_{N}(t_j)-}, V^{(N)}_{\lambda_{N}(t_j)}]$ follows $[
V^{(N)}_{\lambda_{N} (t_i)-}, V^{(N)}_{\lambda_{N}
  (t_i)}]$
for $m_{k+1}>N\ge m_k$. For these $N$-s we define
$
\Lambda_N(T_{t_i-}):= V^{(N)}_{ \lambda_{N} (t_i)-},
$
and 
$
 \Lambda_N(T_{t_i}) = V^{(N)}_{ \lambda_{N} (t_i)}
$
and elsewhere $\Lambda_N(t)$ is defined by a linear
interpolation. 
It is obvious from the construction that $\Lambda_N(t)$ converges uniformly to
$t$ on $[0,T_L]$, as $N\to+\infty$. Combining this with part i) of Lemma \ref{lm3.5} we obtain also that $\lim_{N\to+\infty}|\sigma_N(v)-s(v)|=0$ uniformly on $[0,T_L]$.

Since $\Lambda_N(\cdot)$ dilates  each $[
T_{t_i-}, T_{t_i}]$ onto $[
V^{(N)}_{\lambda_{N} (t_i)-}, V^{(N)}_{\lambda_{N}
  (t_i)}]$ with the scale greater than $d_{k+1}/d_k$  for any $t\in A_{d_k}$ we have
$\Lambda_N(t)\in A_{d_{k+1}}$ and, thanks to Lemma \ref{lm3.5},
we have then for all $N\ge m_{k+1}(\ge k+1)$
\begin{equation}
\label{031802}
S_{\sigma_N(t)} =S_{s\circ
  \Lambda_N (t)} =S_{ s(t)}=S_{t_i},
\end{equation}
see \eqref{012809} for the definition of $\sigma_N(t)$.
The last equality is a consequence of the fact that both
$\Lambda_N(t)$ and $t$ belong to the same $[T_{t_i-}, T_{t_i}]$ for some $i=1,\ldots,\ell_k$ and $s(\Lambda_N(t))=s(t)$. On the other hand we have $u_N\circ \Lambda_N (t)=\lambda_{N} (t_i)$.
Suppose now that $t'\in[T_{t_i-}, T_{t_i}]$ and $t\in[T_{t_i-}, T_{t_i}]\cap A_{d_k}$.
We have then
$
S_{s(t')}=S_{ s(t)}=S_{t_i}
$ and, since $\Lambda_N(t),\Lambda_N(t')\in [
V^{(N)}_{\lambda_{N} (t_i)-}, V^{(N)}_{\lambda_{N}
  (t_i)}]$ we have $u_N(\Lambda_N(t))=u_N(\Lambda_N(t'))=\lambda_N(t_i)$.
  This implies that
  $$
S_{\sigma_N(t')} =S_{\sigma_N(t)} \stackrel{\eqref{031802}}{=}S_{ s(t)}=S_{ s(t')}.
$$
We have shown that
$
\lim_{N\to+\infty}S_{\sigma_N(t)} =S_{ s(t)}
$
uniformly on $s^{-1}({\cal S}_{k})$ for each $k$.

 The statement on the uniform convergence
on the entire $[0, T_L]$ follows from Lemma \ref{uniform}.
Indeed, suppose that $n_k$ is so large that
\begin{equation}
\label{051802}
\sup_{v\in[0,T_L]}(|\sigma_N(v)-s(v)|+|\Lambda_N(v)-v|)<d_k\wedge (\kappa_k/2)\quad\mbox{ for }N\ge n_k.
\end{equation}
Assume also that $t\not\in s^{-1}({\cal S}_{k})$. We claim 
that:
\begin{equation}
\label{061802}
\mbox{ no element from ${\cal S}_{k}$  lies between
$\sigma_N(t)$ and 
$s(t)$ for $N\ge n_k$.}
\end{equation}
 Indeed, suppose that $t_i\in {\cal S}_{k}=s(A_{d_k})$ and
$s(t)<t_i<\sigma_N(t)$.
Then for any $v\in A_{d_k}\cap [T_{t_i-},T_{t_i}]$ we have $v-t>d_k$ and
$$
\lambda_N^{-1}\circ u_N(v)=s(v)=t_i<\sigma_N(t)=\lambda_N^{-1}\circ u_N\circ \Lambda_N (t)
$$ 
and if $v=\Lambda_N(v')$ for some $v'$ then we have to have $v'<t$.
This however implies
$\Lambda_N(v')-v'>v-t>d_k$, which is impossible in light of \eqref{051802}.
The case when $s(t)>t_i>\sigma_N(t)$ can be dealt with similarly.

Suppose that $\varepsilon,\varrho>0$ are arbitrary and $\kappa_k, d_k>0$ are
sufficiently small so that
for a certain $0=s_0<\ldots <s_K=L$ that are $\kappa_k$-sparse we have $ \omega_{S}'(\kappa_k;L) <\ep$ and
$$
\sup_{s_{i-1} \leq x,y <s_i} |S_y-S_x|\le \omega_{S}'(\kappa_k;L) +\varrho.
$$
Recall here that  
$$
\omega_{S}'(\delta;L):=\inf_{t_i}\sup[|S_t-S_{t'}|,\,t_i\le t,t'<t_{i+1},\,i=0,\ldots,N-1], 
$$
where the infimum extends over all partitions $0=t_0<t_1<\ldots<t_N=L$ that are $\delta$-sparse, i.e. such that $\delta<t_{i+1}-t_i$ for all $i=0,\ldots,N-1$, see pp. 109-110 of \cite{billingsley}.

Then, for $n_k$ sufficiently large, so
that \eqref{051802}, thus also claim \eqref{061802}, hold for any $N\ge n_k$
between $ \sigma_N(t)$ and
$s(t)$ there can be at most one $s_i$. Indeed in the case when there were at least two such $s_i$-s we could estimate $|\sigma_N(t)-s(t)|\ge \kappa_k$, which would clearly contradict \eqref{051802}.

  If there is no $s_i$ lying between $\sigma_N(t)$ and $s(t)$  we estimate
$|S_{\sigma_N (t)} -S_{
  s(t)}|\le \om'_{S}(\kappa_k;L)+\varrho$. If, on the other hand, there is such a $s_i$ then according to \eqref{061802} it cannot belong to $ s(A_{d_k})$ and we can estimate
%
$$
|S_{\sigma_N (t)} -S_{
  s(t)}|\le 2[\om'_{S}(\kappa_k;L)+\varrho]+\max_{u\not\in s(A_{d_k})}\Delta S_u<2(\varrho+\varepsilon)+\max_{u\not\in s(A_{d_k})}\Delta S_u.
$$
Summarizing, we have shown that
$$
\limsup_{N\to+\infty}\sup_{t\in[0,T_L]}|S_{\sigma_N (t)} -S_{
  s(t)}|\le 2(\varrho+\varepsilon)+\max_{u\not\in s(A_{d_k})}\Delta S_u
$$
and \eqref{020505} follows.

{\bf Step 2.}  Let $\si_N^*(t):=\lambda_N^{-1}\circ u_N^*\circ \Lambda_N(t)$, where, as we recall $u_N^*(t)=u_N(t)-1/N^{\al}$. Choose an arbitrary  $t\in s^{-1}({\cal S}_k)$ and  suppose that  $t\in [
T_{t_i-}, T_{t_i}]\cap  A_{d_k}$ for some $i=1,\ldots,\ell_k$. Then, as we know, $\Lambda_N(t)\in A_{d_{k+1}}\cap [
V^{(N)}_{\lambda_{N} (t_i)-}, V^{(N)}_{\lambda_{N}
  (t_i)}]$. Thus $s(t)=t_i$ and $u_N^*\circ \Lambda_N(t)= \lambda_{N} (t_i)-1/N^{\alpha}$.
  We conclude therefore that
  $$
  S_{\si_N^*(t)}=S_{s(t)-c_N}
  $$
  for some $c_N>0$ such that $c_N\to0$, as $N\to+\infty$. As a result we get
   \begin{equation}
   \label{022909}
\lim_{N\to+\infty}  S_{\si_N^*(t)}=\zeta^-_t,
  \end{equation}
  fo all $t\in\bigcup_{k\ge1}s^{-1}({\cal S}_k)$. On the other hand, if  $s(t)\not\in \bigcup_{k\ge1}{\cal S}_k$ then, according to Lemma \ref{uniform}, we have $\Delta S_{s(t)}=0$ and, thanks to \eqref{020505} we conclude that
  $$
 \lim_{N\to+\infty} S_{\si_N^*(t)}=\lim_{N\to+\infty} S_{\si_N(t)}=\zeta_t=\zeta_t^-
  $$
  We have shown therefore that \eqref{022909} holds for all $t\ge0$ pointwise and the limiting function is c\`adl\`ag. To finish the proof it suffices only to observe that, thanks to  \eqref{020505}, the sequence of c\`adl\`ag functions $\{S_{\si_N^*(t)},\,t\ge0\}$ has to converge in the $J_1$ topology on $D[0,L]$. Its limit has to coincide with the pointwise limit, because the set of discontinuity points of a c\`adl\`ag function is at most countable, see Corollary 12.2.1 of \cite{whitt}, p. 473. This ends the proof of
  the theorem.
\qed

\section{An application to a jump process arising in a quantum
  transport problem}

\label{quantum}
We illustrate the results obtained in the previous sections with an application
to a jump process that arises in quantum
mechanical transport, see \cite{clark-etal}. Recall that
the one dimensional torus $\bbT$ is an interval $[-\pi,\pi]$ with the
endpoints identified.
Suppose that  $\{K_t,\,t\ge0\}$ is a jump process on $\bbT$ whose generator is given by
\begin{equation}
\label{lind}
Lf(k)=\gamma(k)\int_{\bbT}\hat r(\theta,k)[f(\theta)-f(k)]d\theta
\end{equation}
for $f\in B_b(\bbT)$ - the space of bounded, Borel measurable
functions on $\bbT$. Function $r_0^{-1}\ge\hat r(\theta,k)\ge r_0>0$
is continuous on $\bbT\times \bbT$, even and  {\em doubly stochastic}, i.e.   $\hat r(-\theta,-k)= \hat r(\theta,k)$ and 
$$
\int_{\bbT}\hat r(\theta,k)d\theta=\int_{\bbT}\hat r(k,\theta)d\theta=1,\quad\forall\,k\in\bbT.
$$
 On the other hand, we assume that $\gamma(k)$ is also even and
 strictly positive except for a possible set consisting of two points
 $\{-k_0,k_0\}$. More precisely we have $\gamma(-k)=\gamma(k)$  and
 there is $k_0\in\bbT$ such that  $\gamma(k_0)=0$ and
 $\inf_{|k-k_0|\ge \delta}\gamma(k)>0$ for any $\delta>0$. We suppose
 furthermore that $\gamma(k)\le t_*^{-1}$ for some $t_*>0$
 and $\int_{\bbT}\gamma^{-1}(k)dk=+\infty$.
This kind of processes appears while considering the transport of particles in quantum systems, see e.g. Section 4.3 of \cite{clark-etal}. 
It is easy to see 
that $m_*(dk)=\gamma^{-1}(k)m_1(dk)$ is an infinite, reversible, invariant measure for the process.  Here $m_1(dk)=dk/(2\pi)$ denotes the normalized Lebesgue measure on the torus. 
 Indeed, for any $f\in B_b(\bbT)$
\begin{eqnarray*}
&&
\int_{\bbT} Lf(k)m_*(dk)=\frac{1}{2\pi}\int_{\bbT} \int_{\bbT}\hat r(\theta,k)[f(\theta)-f(k)]d\theta dk\\
&&
=\frac{1}{2\pi}\int_{\bbT}f(\theta) d\theta\int_{\bbT} \hat r(\theta,k)dk-\frac{1}{2\pi}\int_{\bbT} \int_{\bbT}\hat r(\theta,k)f(k)d\theta dk=0.
\end{eqnarray*}
The process $\{K_t,\,t\ge0\}$ can be realized using a Markov chain and a renewal process that corresponds to the jump times. Consider a skeleton Markov chain $\{(X_n,\rho_n),\,n\ge0\}$, defined on $\bbT_{k_0}\times (0,+\infty)$, where $\bbT_{k_0}:=\bbT\setminus\{-k_0,k_0\}$, $\{\rho_n,\,n\ge0\}$ is an i.i.d. sequence of exponentially distributed random variables with intensity $1$ and  $\{X_n,\,n\ge0\}$ is an independent  Markov chain   with the state space $\bbT_{k_0}$,  whose  transition probability equals $\hat r(\theta,k)d\theta$.
Let $\tau(k,\rho):=\gamma^{-1}(k)\rho$, $t_0:=0$ and
$t_n:=\sum_{k=0}^{n-1}\tau(X_k,\rho_k)$, $n\ge1$. We let $K_t:=X_n$ for $t\in [t_n,t_{n+1})$.

\subsection{Harris recurrence property}
Our first result concerns the recurrence property of $\{K_t,\,t\ge0\}$.
\begin{prop}
\label{prop-harris}
Suppose that $h\in(0,t_*)$.  Consider an embedded Markov chain $\{K_{nh},\,n\ge0\}$. It is Harris recurrent w.r.t.  measure $m_1$, i.e. for any Borel subset $B$ with $m_1[B]>0$ we have
\begin{equation}
\label{0.1}
\bbP[\exists\,n\ge0:\, K_{nh}\in B]=1.
\end{equation}
\end{prop}
\proof
To simplify the notation let $h=1$.
Our hypotheses on the skeleton chain guarantee that
\begin{equation}
\label{010209}
\bbP[D]=1,
\end{equation}
where $D=[X_n\in B,\,i.o.]$.
Let $A_n:=[X_n\in B,\,t_{n+1}-t_n\ge 2]$.
To see that \eqref{0.1} holds it suffices only to prove that
\begin{equation}
\label{0.3}
\bbP[C]=1,
\end{equation}
where $C:=\bigcup_{n\ge0}A_n$.
Note that $1_{C^c}(\om)\le f(\om)$, where
$$
f(\om):=\prod_{n\ge0}
\left[1_{B^c}(X_n)+1_{B}(X_n)1_{[\tau(X_n,\rho_n)< 2]}\right].
$$
However, 
\begin{eqnarray*}
&&
\bbE f=\bbE\left\{\prod_{n\ge0}
\left[1_{B^c}(X_n)+1_{B}(X_n)(1-e^{-2\gamma(X_n)})\right]\right\}
\\
&&
\le \bbE\left\{\prod_{n\ge0}
\left[1_{B^c}(X_n)+1_{B}(X_n) (1-e^{-2/t^*})\right]\right\}
\\
&&
=\bbE\left\{\prod_{n\ge0}
\left[1_{B^c}(X_n)+1_{B}(X_n) (1-e^{-2/t^*})\right],\,D\right\}
=0
\end{eqnarray*}
and \eqref{0.3} follows in light of \eqref{010209}.
\qed

As an immediate corollary from the above proposition and Theorem 1,
p.~116 of \cite{harris} we obtain that $m_*$ is the unique $\si$-finite invariant measure under the process that is absolutely continuous w.r.t. $m_1$.

 Denote by $\{P_t,\,t\ge0\}$ the transition semigroup of the process $\{K_t,\,t\ge0\}$. 
It satisfies the following integral equation
\begin{equation}
\label{int}
P_tf(k)=e^{-t\gamma(k)}f(k)+\gamma(k)\int_0^te^{-s\gamma(k)}ds\int_{\bbT}\hat r(k',k)P_{t-s}f(k')dk'.
\end{equation}
For any $N\ge 1$ and $T>0$ denote by 
$$
\Delta_N(T):=[(s_0,\ldots,s_{N-1}):s_i\ge0,\,i=0,\ldots,N-1,\,\sum_{i=0}^{N-1}s_i\le T].
$$
Iterating equation \eqref{int} we can easily show that
\begin{eqnarray*}
&&
P_tf(k)=e^{-t\gamma(k)}f(k)+\sum_{N=1}^{+\infty}\gamma(k)\mathop{\int\ldots\int}_{\Delta_N(t)}\mathop{\int\ldots\int}_{(\bbT)^{N}}e^{-t\gamma(k_N)}
\\
&&
\times\prod_{i=1}^N\left\{\gamma(k_i)e^{-s_i(\gamma(k_i)-\gamma(k_N))}\hat r(k_i,k_{i-1})\right\}f(k_N)ds^{(N)}dk^{(N)}.
\nonumber
\end{eqnarray*}
Here $k_0:=k$, $ds^{(N)}:=ds_0\ldots ds_{N-1}$, $dk^{(N)}:=dk_1\ldots dk_N$. 
The  component of transition probability that is absolutely continuous w.r.t. $m_1$
equals therefore
\begin{eqnarray*}
\!\!\!\!&&
p_t(k,k')=\sum_{N=1}^{+\infty}\gamma(k)\!\!\mathop{\int\!\!\ldots\!\!\int}_{\Delta_N(t)}\!\!\mathop{\int\!\!\ldots\!\!\int}_{(\bbT)^{N-1}}e^{-t\gamma(k')}
\prod_{i=1}^N\left\{\gamma(k_i)e^{-s_i(\gamma(k_i)-\gamma(k'))}\hat r(k_i,k_{i-1})\right\}ds^{(N)}dk^{(N-1)}.
\nonumber
\end{eqnarray*}
Here $k':=k_N$.Thus, for
every $h>0$, $C\in{\cal B}(\bbT_{k_0})$ with dist$(C,\{-k_0,k_0\})>0$ and $m_1(C)>0$ we have
$
\inf_{k,k'\in C}p_{h}(k,k')>0.
$
The transition probability function of any embedded chain
$\{K_{nh},\,n\ge0\}$ is therefore aperiodic in the sense of \cite{jamison-orey}.
Suppose that $f_0=d\nu_0/dm_*\in L^2(m_*)$ is a density.
Thanks to reversibility of $m_*$ we obtain that $\nu_0P_t$ is absolutely continuous w.r.t. $m_1$ and its density equals
$$
f_t=\frac{d\nu_0P_t}{dm_*}=P_tf_0,\quad\forall\, t\ge0.
$$

\subsection{Mixing property of the process} 
\begin{thm}
\label{thm1}
Suppose that the initial law $\nu_0$ is absolutely continuous w.r.t. the Lebesgue measure $m_1$.
Then, $\nu_0P_t$ converge weakly, as $t\to+\infty$, to the measure
$\mu_*:=1/2(\delta_{k_0}+\delta_{-k_0})$. In addition, the process is
completely mixing, i.e. if $\nu_0$, $\nu_0'$ are two absolutely continuous initial laws then 
\begin{equation}
\label{0.11}
\lim_{t\to+\infty}\|\nu_0 P_t-\nu_0' P_t\|_{TV}=0.
\end{equation}
\end{thm}
To prove the above result we first show the following.
\begin{prop}
\label{prop1}
For any compact set $K\subset \bbT_{k_0}$ and a measure $\nu_0$
as 
in Theorem \ref{thm1} we have
\begin{equation}
\label{sweeping}
\lim_{t\to+\infty}\nu_0 P_t[K]=0.
\end{equation}
\end{prop}
\proof
Using a density argument it suffices to show \eqref{sweeping} 
for measures $\nu_0$ whose density belongs to $L^2(m_*)$. Thanks to strong
continuity of semigroup $\{P_t,\,t\ge0\}$ in $L^1(m_*)$ in order to prove
\eqref{sweeping} it suffices only to show
that for any $h>0$
\begin{equation}
\label{sweeping-1}
\lim_{n\to+\infty}\nu_0 P_{nh}[K]=0.
\end{equation}
From the Harris recurrence property, see Proposition \ref{prop-harris}, we know  that for any set $A\subset \bbT_{k_0}$ with $m_*[A]>0$ we have $P_h1_A(x)>0$, $m_*$-a.e. hence from
\cite{foguel}, pp. 85-102, we have \eqref{sweeping-1} for any $K$ such that $+\infty>m_*[K]>0$,
cf. Theorem C, p. 91 of ibid.
\qed

{\bf The proof of Theorem \ref{thm1}.}
Let $h\in(0,t_*)$.
Define by ${\cal C}:=\bigcap_{n\ge0}{\cal C}_n$, where ${\cal C}_n$ is the smallest $\si$-algebra generated by $\{K_{mh},\,m\ge n\}$.
According to Theorem 1, p. 45 of \cite{jamison-orey} the tail $\si$-algebra of the chain that is Harris recurrent and aperiodic has to be trivial. Therefore, according to Lemma 3, p. 43 of ibid. \eqref{0.11} follows.

Observe that if $u(k)$ is a density w.r.t. $m_*(dk)=\gamma^{-1}(k)dk$ such that $u(-k)=u(k)$ then 
\begin{equation}
\label{even}
uP_t(-k)=uP_t(k).
\end{equation}    This follows from the fact that
$v_t(k):=uP_t(-k)$ satisfies equation 
\begin{eqnarray*}
&&\dfrac{dv_t}{dt}(k)
=
v_tL(k),\quad
v_0(k)=u(k).
\end{eqnarray*}
Since $uP_t(k)$ satisfies the same equation from the uniqueness of solutions we obtain $uP_t(k)=v_t(k)=uP_t(-k)$.
Let $\nu_0(dk):=u(k)m_*(dk)$. Combining \eqref{even}  with Proposition
\ref{prop1} and \eqref{0.11} we conclude that $\nu_0P_t\Rightarrow
(1/2)(\delta_{-k_0}+\delta_{k_0})$, as $t\to+\infty$. From the
(already shown) complete mixing property we conclude in particular
that for any initial distribution $\mu_0$, absolutely continuous
w.r.t. $m_1$, we have  $\mu_0P_t\Rightarrow (1/2)[\delta_{-k_0}+\delta_{k_0}]$, weakly over $C(\bbT)$, as $t\to+\infty$.\qed

\subsection{Limit theorems for additive functionals of the process}

In this section we shall be concerned with the convergence of the laws of $N^{-\ga}\int_0^{Nt}\Psi(K_s)ds$ for an appropriate $\ga>0$ and $\Psi(k)$.
Suppose that $\gamma(k)\sim c_*|k-k_0|^{\kappa}$, when $|k-k_0|\ll 1$,
for some $\kappa>1$ and $c_*>0$. Then the law of
$\tau(k,\rho)=\gamma^{-1}(k)\rho$ under $m_1\otimes \lambda$, where $\lambda(d\rho):=e^{-\rho}d\rho$, belongs to the domain of
attraction of a stable subordinator with index
$\alpha=1/\kappa$ so Condition  \ref{tails} holds. One can easily verify that the other assumptions  about the Markov chain made in Section \ref{markov-chain} hold for  $\{(X_n,\rho_n),\,n\ge0\}$ as well.

 Note that 
 $
 \int_0^{t}\Psi(K_s)ds=\widehat W(t),
 $ where $\widehat W(t)$ is the linear interpolation of a CTRW corresponding to renewal times $\{t_N,\,N\ge0\}$ and the partial sums $\{S_N,\,N\ge0\}$ formed for 
 $
 V(k,\rho):=\Psi(k)\tau(k,\rho).
 $
When $V\in L^2(m_1\otimes \lambda)$ and 
\begin{equation}
\label{centering}
\int_{\bbT}\int_0^{+\infty}V(k,\rho)dk\la(d\rho)=0
\end{equation}
  from Theorem \ref{dependent} it follows that:
\begin{corollary}
\label{cor022101}
The processes $\{Y^{(N)}_t:=N^{-\al/2}\int_0^{Nt}\Psi(K_s)ds,\,t\ge0\}$ converge over $C[0,+\infty)$ converge, as $N\to+\infty$,  to the law of the Mittag-Leffler process that corresponds to an $\al$-stable subordinator.
\end{corollary}

Assume also that the law of $V(k,\rho)$ under $m_1\otimes \la$ belongs to the
domain of attraction of a $\beta$-stable law.  Denote by $s_N(t)$ the right inverse of $\{T^{(N)}_t,\,t\ge0\}$
\begin{corollary}
\label{cor022101b}
Suppose that $\beta\not=1$ and for $\beta\in(1,2)$ condition \eqref{centering} holds.
In addition, assume that $\tau(k,\rho)$ and $V(k,\rho)$ satisfy \eqref{sing-3}. 
Then,  
$\{Y^{(N)}_t:=N^{-\al/\beta}\int_0^{Nt}\Psi(K_s)ds,\,t\ge0\}$ converge in law over ${\cal D}$ with the $M_1$ topology 
to   
$\{\zeta_t:=S_{s(t)},\,t\ge0\}$, where $\{S_t,\,t\ge0\}$ is a
$\beta$-stable process and $s(t)$ is the right inverse of an independent,
$\al$-stable subordinator. When $\beta=1$ the theorem still holds
for  $\{Y^{(N)}_t-c_N s_N(t),\,t\ge0\}$, where 
$
 c_N:=\int_{[|V(k,\rho)|\le N]}V(k,\rho)dk\la(d\rho).
 $
\end{corollary}

{\bf Remark.} The jump process considered in Section 4.3 of \cite{clark-etal} 
has the generator given by $Lf(k)=c\cos^2k\int_{\bbT}\hat
r(k'-k)[f(k')-f(k)]dk'$ for some constant $c>0$ and a  density
function $\hat r(k)$ satisfying $r_*\le \hat r(k)\le r_*^{-1}$ for
some $r_*\in(0,1)$. Theorem \ref{thm1} implies that $\nu_0
P_t\Rightarrow 1/2(\delta_{-\pi/2}+\delta_{\pi/2})$, as $t\to+\infty$
for any initial measure $\nu_0$ absolutely continuous w.r.t. 
the Lebesgue measure. This answers in the affirmative the conjecture made in the discussed paper. 
For an observable $\Psi$ such that $V=\Psi\tau\in L^2(m_1\otimes \la)$ the functionals $\{N^{-1/4}\int_0^{Nt}V(K_s)ds,\,t\ge0\}$ converge in law to a Mittag-Leffler process.
Note that when $\Psi(-k)=-\Psi(k)$ belongs to the normal domain of attraction of a Cauchy law and it is not singular at $\pi/2$ then $\{N^{-1/2}\int_0^{Nt}\Psi(K_s)ds,\,t\ge0\}$ converge in law, as $N\to+\infty$, to a process that has the same scaling properties as a Brownian motion but is not Markovian (recall that we call it a fake diffusion).

\begin{appendix}

\section{The proofs of Theorems \ref{joint-conv-2} and \ref{joint-conv-1}}

\label{appA}

Theorems in question can be
concluded from the result we formulate below. It is essentially a two dimensional version of Theorem 4.1 p. 840 of  \cite{durrett}.
Since its  proof is a simple modification of an argument presented ibid.   we omit its presentation here.
Also, for simplicity we only consider the case when $K_N=N$.

Before formulating the result we introduce some notation. Suppose that $\{Z_{n,N},\,n\ge0,\,N\ge1\}$ is an array of $\bbR^2$-valued random vectors on $(\Om,{\cal F},\bbP)$ and $\{{\cal G}_{n,N},n\ge{-1},\,N\ge1\}$ is an array of sub $\si$-algebras of ${\cal F}$ such that ${\cal G}_{n-1,N}\subset {\cal G}_{n,N}$, $N\ge1$, $n\ge0$.
For a fixed $\Delta=(\Delta_1,\Delta_2)\in(0,+\infty)^2$ define $
Z^{(\Delta,N)}_t:=\sum_{n=0}^{[Nt]-1}Z_{n,N}^{(\Delta)}
$, where
$$
Z_{n,N}^{(\Delta)}:=Z_{n,N}1_{[|Z_{n,N}^{(i)}|\le \Delta_i,i=1,2]}-\bbE\left[Z_{n,N}1_{[|Z_{n,N}^{(i)}|\le \Delta_i,i=1,2]}\Big|{\cal G}_{n-1,N}\right]
$$
and
$$
\tilde Z^{(\Delta,N)}_t:=\sum_{n=0}^{[Nt]-1}\left[Z_{n,N}-\bbE\left[Z_{n,N}1_{[|Z_{n,N}^{(i)}|\le \Delta_i,i=1,2]}\Big|{\cal G}_{n-1,N}\right]\right].
$$
Let $|(x_1,x_2)|_{\infty}:=\max\{|x_1|,|x_2|\}$. 
 Suppose  that $\nu_0$ is a measure on $\mathbb R^2_*:=\mathbb R^2\setminus\{(0,0)\}$ such that $\int_{\mathbb R^2_*}(|x|^2\wedge 1)\nu_0(dx)<+\infty$. To simplify the statement we assume that $\nu_0$ is absolutely continuous w.r.t. the two dimensional Lebesgue measure. 
\begin{thm} (See \cite{durrett})
\label{prop010805}
Assume  that for any $g\in C_b^\infty(\mathbb R^2_*)$ such that $0\not\in${\em supp }$g$ we have:
\begin{equation}
  \label{eq:16b}
  \lim_{N\to+\infty}\mathbb E\left|\sum_{n=0}^{N}\mathbb E\left[g\left(
        Z_{n,N}\right)\Big|{\cal G}_{n-1,N}\right]-
    \int_{\mathbb R^2_*}g(\la_1,\la_2)\nu_0(d\la_1,d\la_2) \right|=0,
\end{equation}
 \begin{equation}
   \label{eq:16bc}
   \lim_{N\to+\infty}\sum_{n=0}^N\mathbb E\left\{\mathbb E\left[g\left(
        Z_{n,N}\right)\Big|{\cal G}_{n-1,N}\right]\right\}^2=0
  \end{equation}
  and
  \begin{equation}
  \label{eq:16bcde}
  \lim_{|\Delta|_{\infty}\to 0+}\limsup_{N\to+\infty}\mathbb E\left[\sup_{t\in[0,T]}|  Z^{(\Delta,N)}_t|^2_{\infty}\right]=0
\end{equation}
 for any $T>0$.
Then, for each $\Delta=(\Delta_1,\Delta_2)\in(0,+\infty)^2$ the processes $\{\tilde Z_t^{(\Delta, N)},\,t\ge0\}$ converge in law over ${\cal D}_2$, with the $J_1$-topology, to the L\'evy process $\{Z_t=(Z_t^{(1)},Z_t^{(2)}),\,t\ge0\}$ such that
$
\bbE e^{i\xi_1 Z_t^{(1)}+i\xi_2 Z_t^{(2)}}=e^{t\psi(\xi_1,\xi_2)}
$
and
\begin{eqnarray}
\label{032602}
\psi(\xi_1,\xi_2):=\int_{\bigcup_{i=1}^2[|\la_i|\ge \Delta_i]}(e^{i\xi_1 \la_1+i\xi_2 \la_2}
-1)\nu_0(d\la_1,d\la_2)\\
+\int_{[|\la_i|\le \Delta_i,\,i=1,2]}[e^{i\xi_1 \la_1+i\xi_2 \la_2}
-1-i(\la_1\xi_1+\la_2\xi_2)]\nu_0(d\la_1,d\la_2).
\end{eqnarray}
\end{thm}

To use the above  theorem we transform slightly the process $\{X_t^{(N)}:=(S_t^{(N)},T_t^{(N)}),\,t\ge0\}$. 
In order to simplify the notation we  consider only the case $\beta\in(1,2)$, the
consideration in the case $\beta\in(0,1]$ can be
done similarly (in fact it is even  simpler then).  
Suppose $\ga\in(1,\beta)$. Thanks to the spectral gap condition \eqref{eq:8}  we can find a  unique zero mean solution $\chi$ in $L^{\ga}(\pi)$ of 
\begin{equation}
\label{eq:14}
(I-P)\chi=V.
\end{equation}
For a fixed $M>0$ we let $M_N:=M N^{1/\al}$ and 
\begin{equation}
\label{010103}
\tau^{(N)}(x):=\tau(x)1_{[\tau(x)<M_N]}.
\end{equation} 
 We let also $Z_{n,N}(M):=(Z_{n,N}^{(1)},Z_{n,N}^{(2)}(M))$, where
\begin{eqnarray}
\label{010805}
&&
Z_{0,N}^{(1)}:=0,\quad Z_{n,N}^{(1)}:=\frac{1}{N^{1/\beta}}R_0(X_n,X_{n-1}),\,n\ge1.
\end{eqnarray}
Here $R_0(x,y):=\chi(x)-P\chi(y)$. 
In addition,
$ Z_{n,N}^{(2)}(M):=N^{-1/\al}
\tau^{(N)}(X_n),\,n\ge0$. 
For $N,n\ge0$ we let ${\cal G}_{n,N}$ be the $\si$-algebra generated by $X_0,\ldots,X_n$.
By convention we let ${\cal G}_{-1,N}$ be the trivial $\si$-algebra.

Define a process
\begin{equation}
\label{122702}
Z_t^{(N)}(M)=(Z_{t,N}^{(1)},Z_{t,N}^{(2)}(M)):=\sum_{n=0}^{[Nt]-1}Z_{n,N}(M),\quad t\ge0.
\end{equation}
 A simple calculation shows that
\begin{eqnarray*}
&&X_t^{(N)}-Z_t^{(N)}(M)=\left(N^{-1/\beta}[\chi(X_0)-P\chi(X_{[Nt]-1})],N^{-1/\al}\sum_{n=0}^{[Nt]-1}\tau(X_n)1_{[\tau(X_n)>M_N]}\right).
\end{eqnarray*}
Hence, due to \eqref{062202} (proved below) and \eqref{eq:tails}, for any $T,\sigma>0$ we have
\begin{equation}
\label{012702}
\lim_{M\to+\infty}\limsup_{N\to+\infty}\bbP[\sup_{t\in[0,T]}|X_t^{(N)}-Z_t^{(N)}(M)|_{\infty}\ge \sigma]=0.
\end{equation}
To prove Theorems \ref{joint-conv-2} and \ref{joint-conv-1} it   suffices therefore to show the convergence of processes $\{Z_t^{(N)}(M),\,t\ge0\}$, $N\ge1$ for a fixed $M>0$ to a L\'evy process $\{Z_t(M),\,t\ge0\}$ whose exponent equals
\begin{equation}
\label{052602}
\psi(\xi_1,\xi_2)=
\int_{\mathbb R^2}[e^{i\xi_1 \la_1+i\xi_2 \la_2}
-1-i(\la_1\xi_1+\la_2\xi_2)]1_{[|\la_2|\le M]}\nu_*(d\la_1,d\la_2).
\end{equation}
Since now on we drop $M$ from our subsequent notation of stochastic processes writing $Z_{n,N}:=Z_{n,N}(M)$, $Z_t^{(N)}:=Z_t^{(N)}(M)$ and $Z_t:=Z_t(M)$.

Below we verify that  processes $\{Z_t^{(N)},\,t\ge0\}$, $N\ge1$  satisfy the assumptions of Theorem
\ref{prop010805} with measure  $\nu_0$   equal $\nu_*$ from either Theorem \ref{joint-conv-2}, or \ref{joint-conv-1} correspondingly. Accepting this claim for a moment we show how to reach the conclusions of the aforementioned theorems. Let $\Delta_2:=M$. 
From  Theorem
\ref{prop010805}  we deduce that 
\begin{equation}
\label{042602}
\tilde Z_t^{(N)}(\Delta_1):=Z_t^{(N)}-\sum_{n=0}^{[Nt]-1}\bbE\left[Z_{n,N}1_{[|Z_{n,N}^{(1)}|\le \Delta_1]}1_{[|Z_{n,N}^{(2)}|\le \Delta_2]}\Big|{\cal G}_{n-1,N}\right],\,t\ge0\}
\end{equation} converge in law over ${\cal D}_2$ to 
a L\'evy process $\{Z_t,\,t\ge0\}$ whose exponent equal $\psi(\xi_1,\xi_2)$ given by \eqref{032602} (with $\Delta_2=M$). 
Since $\bbE\left[Z_{n,N}^{(1)}\Big|{\cal G}_{n-1,N}\right]=0$ and $|Z_{n,N}^{(2)}|\le M$ a.s. we get
\begin{eqnarray}
\label{012602}
\sum_{n=0}^{[Nt]-1}\bbE\left[ Z_{n,N}^{(1)}1_{[|Z_{n,N}^{(i)}|\le \Delta_i,\,i=1,2]}\Big|{\cal G}_{n-1,N}\right]=
-\sum_{n=0}^{[Nt]-1}\bbE\left[ Z_{n,N}^{(1)}1_{[|Z_{n,N}^{(1)}|\ge \Delta_1]}\Big|{\cal G}_{n-1,N}\right].
\end{eqnarray}
This together with  \eqref{062202}, \eqref{eq:tails:1b} (proved below) and hypothesis \eqref{eq:tails} imply  that for any $T>0$ we have
\begin{eqnarray}
\label{012602b}
\lim_{\Delta_1\to+\infty}\limsup_{N\to+\infty}\bbE\left|\sup_{t\in[0,T]}\sum_{n=0}^{[Nt]-1}\bbE\left[Z_{n,N}1_{[|Z_{n,N}^{(1)}|\le \Delta_1]}1_{[|Z_{n,N}^{(2)}|\le M]}\Big|{\cal G}_{n-1,N}\right]\right|=0.
\end{eqnarray}
Letting $\Delta_1\to+\infty$ we deduce that $\{Z_t^{(N)},\,t\ge0\}$ converge in law over ${\cal D}_2$ to a L\'evy process with L\'evy exponent \eqref{052602}.


We start with the following lemma that, among others, allows us to justify the limits in \eqref{012702} and \eqref{012602}.
\begin{lemma}
 \label{lm012302}
 We have 
 \begin{equation}
    \label{eq:tails:1b}
    \begin{split}
      & \pi(\chi\ge\lambda)    =  \frac{c_\beta^+}{\lambda^{\beta}}\left(1+o(1)\right),
      \\
      &\pi(\chi\le-\lambda)    =  \frac{c_\beta^-}{\lambda^{\beta}}\left(1+o(1)\right),\quad \mbox{as $\lambda\to+\infty$. }
    \end{split}
  \end{equation}
In addition,  
\begin{equation}
\label{062202}
\|P|\chi|\|_{L^{2}(\pi)}<+\infty
\end{equation}
and there exists $C>0$ such that
\begin{equation}
\label{022602}
\|P\tau^{(N)}\|_{L^2(\pi)}\le CN^{1/\al-1},\quad\forall\,N\ge1,\,M>0.
\end{equation}
 \end{lemma}
 \proof
 Observe first that
\begin{eqnarray}
\label{032202}
&&\|P|V|\|_{L^2(\pi)}^2\stackrel{\eqref{new-a}}{\le} C\int\left(\int_0^{+\infty} P_a(x,|V|>\la)d\la\right)^2\pi(dx)
\\
&&
= C\int\left(\int p(x,y)|V(y)|\pi(dy)\right)^2\pi(dx)\nonumber\\
&&
\stackrel{\mbox{\tiny Jensen}}{\le} 
C\int\left(\int p^2(x,y)\frac{|V(y)|\pi(dy)}{\|V\|_{L^1(\pi)}}\right)\pi(dx)\|V\|_{L^1(\pi)}^2\le CC(2)\|V\|_{L^1(\pi)}^2.\nonumber
\end{eqnarray}
Since
    $
\chi=(I-P)^{-1}V=\sum_{n\ge0}P^nV
$
from the above we deduce that $P|\chi|\le (I-P)^{-1}P|V|$, thus
 $P|\chi|\in L^2(\pi)$. From this and the Poisson equation
\eqref{eq:14} it follows that $\chi$ satisfies \eqref{eq:tails:1b}.

To show \eqref{022602} we  write
\begin{eqnarray}
\label{042202}
&&\|P\tau^{(N)}\|_{L^2(\pi)}^2    \le
    C(2)\|\tau^{(N)}\|_{L^1(\pi)}^2\le    C(2)\left(\int_{0}^{M_N} \pi[\tau>\la]d\la\right)^2
\\
&&
\le C\left(\int_{0}^{M_N} \frac{ d\la}{1+\la^{\al}}\right)^2\le C' N^{2(1/\al-1)}.\qed\nonumber
    \end{eqnarray}

\subsubsection*{The proof of \eqref{eq:16b}} Since
\begin{eqnarray*}
&&
\lim_{\Delta_1\to+\infty}\limsup_{N\to+\infty}\bbE\left|\sup_{t\in[0,T]}\sum_{n=0}^{[Nt]-1}\bbP\left[|Z_{n,N}^{(1)}|\ge \Delta_1,\mbox{or }\,|Z_{n,N}^{(2)}|\ge M\Big|{\cal G}_{n-1,N}\right]\right|\\
&&
=\lim_{\Delta_1\to+\infty}\limsup_{N\to+\infty}\bbE\left|\sup_{t\in[0,T]}\sum_{n=0}^{[Nt]-1}\bbP\left[|Z_{n,N}^{(1)}|\ge \Delta_1\Big|{\cal G}_{n-1,N}\right]\right|=0
\end{eqnarray*}
it suffices to show \eqref{eq:16b} for $g\in C^{\infty}_0(\bbR^2_*)$.
In that case we can expand $g\left(
Z_{n+1,N}\right)$ using Taylor formula around
$z^{(N)}(X_{n+1})$, where 
$
z^{(N)}(x):=(N^{-1/\beta} V(x),N^{-1/\al}
\tau^{(N)}(x)),
$
and  obtain that
\begin{eqnarray}
  \label{eq:170bb}
&&  \sum_{n=0}^{N-1}\mathbb E\left[g\left(
Z_{n+1,N}\right)\left|\right.{\cal G}_n\right]=\sum_{n=0}^{N-1}\mathbb E\left[g\left(z^{(N)}(X_{n+1})\right)\left|\right.{\cal G}_n\right]
\\
&&
+ \sum_{n=0}^{N-1}\int_0^1\mathbb E\left[R(X_{n+1},X_n)\cdot\nabla g\left(
 z_n^{(N)}(\la)\right)\left|\right.{\cal G}_n\right]d\la .
\nonumber
\end{eqnarray}
Here
$
z_n^{(N)}(\la):=\la R(X_{n+1},X_n)+z^{(N)}(X_{n+1})
$
and
$$
R(x,y):=(N^{-1/\beta}[P\chi(x)-P\chi(y)],-N^{-1/\al}
P\tau^{(N)}(y)).
$$
Denote the first and  second terms on the right hand side of  \eqref{eq:170bb} by $I_N$ abd $I\!I_N$
respectively.
We can write $I_N=I_N^{(1)}+I_N^{(2)}$, where
\begin{eqnarray*}
&&
I_N^{(1)}:=\sum_{n=1}^{N}\int g(z^{(N)}(y))P(X_{n-1},dy)-N\int g\left(z^{(N)}(x)\right)\pi(dx),\\
&&
I_N^{(2)}:=N\int g\left(z^{(N)}(x)\right)\pi(dx).
\end{eqnarray*}
Note that 
\begin{eqnarray}
\label{082801b}
&&
\bbE|I_N^{(1)}|=\bbE\left|\sum_{n=1}^{N}PG_N(X_{n-1})\right|,
\end{eqnarray}
where
$
G_N(x):=g(z^{(N)}(x))-\int g\left(z^{(N)}(y)\right)\pi(dy).
$
 Let $u_N := (I-P)^{-1} PG_N$. By the spectral
gap condition (\ref{eq:9}) we have
\begin{equation*}
\label{082802b}
  \int u_N^2(x) \; \pi(dx) \le \frac 1{1-a^2}\int (PG_N)^2(x) \; \pi(dx).
\end{equation*}
Let $a_*:=$dist$(0,$supp~$g)$. 
Note that
\begin{eqnarray}
\label{082202}
&&
\int (PG_N)^2(x) \; \pi(dx)\le 2\|g\|_{\infty}^2\left[\int \left(\int_{[\tau(y)\ge a^*N^{1/\al}/2]}p(x,y)\pi(dy)\right)^2 \; \pi(dx)\right.\\
&&
\left.+\int \left(\int_{[|V(y)|\ge a^*N^{1/\beta}/2]}p(x,y)\pi(dy)\right)^2 \; \pi(dx)\right]\nonumber\\
&&
\le 2C\|g\|_{\infty}^2N^{-1}
\left[\int \int_{[\tau(y)\ge a^*N^{1/\al}/2]}p^2(x,y)\pi(dy) \; \pi(dx)\right.\nonumber\\
&&
\left.+\int \int_{[|V(y)|\ge a^*N^{1/\beta}/2]}p^2(x,y)\pi(dy) \; \pi(dx)\right]\le 2C\|g\|_{\infty}^2N^{-1}o(1),\nonumber
\end{eqnarray}
as $N\to+\infty$. In addition
we can rewrite 
\begin{equation*}
   \sum_{n=1}^{N} PG_N(X_{n-1}) = u_{N}(X_0) - u_{N}(X_{[Nt]}) + \sum_{n=1}^{[Nt]-1}U_n
 ,
\end{equation*}
where $U_n:=u_{N}(X_n) - Pu_{N}(X_{n-1})$, $n\ge1$ is a stationary sequence of  martingale differences with respect to  $\{{\cal G}_n,\,n\ge0\}$. Consequently,
\begin{equation}
\label{072202}
\bbE|I_N^{(1)}|\le \left\{  \mathbb E\left[ \sum_{n=1}^{N} PG_N(X_{n-1})\right]^2\right\}^{1/2} \le C\left( N  \int (PG_{N})^2(y)
  \; \pi(dy)\right)^{1/2}\to0, 
\end{equation}
by virtue of \eqref{082202}.
To prove that
\begin{equation}
\label{022010b}
\lim_{N\to+\infty}\bbE\left|I_N-\int g(\la_1,\la_2) 1_{[0,M]}(\la_2)\nu_*(d\la_1,d\la_2)\right|=0.
\end{equation}
it suffices to use  the  argument above and the following.
\begin{prop}
\label{lm021701}
Suppose that $g\in C_0^\infty(\bbR^2_*)$.
Then, under the assumptions of either Theorem \ref{joint-conv-2}, or \ref{joint-conv-1},
we have
\begin{eqnarray}
\label{021701b}
&&\lim_{N\to+\infty}I_N^{(2)}
=\int g(\la_1,\la_2) 1_{[0,M]}(\la_2)\nu_*(d\la_1,d\la_2),
\end{eqnarray}
where
$
\nu_*(d\la_1,d\la_2)$
is given by \eqref{021902}, or \eqref{031902} respectively.
\end{prop}
\proof {\em The case \eqref{sing-3} holds.}
Suppose that $\gamma>\kappa_1>\alpha\vee \beta$, where $\gamma$ is the
same as in \eqref{sing-3}, and
\begin{eqnarray*}
&&
A_N:=[\tau^{(N)}\ge (a_*/2)N^{1/\kappa_1} ,\,|V|\ge  (a_*/2) N^{1/\beta}],\\
&&
 B_N:=[\tau^{(N)}\ge  (a_*/2) N^{1/\alpha} ,\,|V|\ge (a_*/2)  N^{1/\kappa_1}].
\end{eqnarray*}
Observe that
$\pi(A_N)= o(1/N)$ and  $\pi(B_N)=o(1/N)$.
To compute the limit in \eqref{021701b} it suffices therefore to compute $\lim_{N\to+\infty}{\cal K}_N^{(i)}$, $i=1,2$, where
\begin{eqnarray*}
&&{\cal K}_N^{(i)}:=N\int g\left(z^{(N)}(x)\right)1_{C_N^{(i)}}\pi(dx),\quad i=1,2,
\end{eqnarray*}
where
\begin{eqnarray*}
&&
C_N^{(1)}:=[\tau^{(N)}\le (a_*/2)N^{1/\kappa_1},\,|V|\ge  (a_*/2) N^{1/\beta}],\\
&&
C_N^{(2)}:=[\tau^{(N)}\ge  (a_*/2) N^{1/\alpha} ,\,|V|\le (a_*/2)  N^{1/\kappa_1}].
\end{eqnarray*}
Up to a term of order $o(1)$ we have ${\cal K}_N^{(i)}=\tilde {\cal K}_N^{(i)}$, where
$
\tilde {\cal K}_N^{(i)}:=N\int g\left(\tilde z^{(N)}_i(x)\right)\pi(dx),
$
and 
$
\tilde z^{(N)}_1(x):=(N^{-1/\beta}V(x),0),
$
$
\tilde z^{(N)}_2(x):=(0,N^{-1/\al}\tau^{(N)}(x)).
$
We can write
\begin{eqnarray*}
&&
\lim_{N\to+\infty}\tilde {\cal K}_N^{(1)}=\lim_{N\to+\infty}N^{1-1/\beta}\int_0^{+\infty} \partial_1 g\left(N^{-1/\beta}\la,0\right)\pi(V>\la)d\la\\
&&
+\lim_{N\to+\infty}N^{1-1/\beta}\int_0^{+\infty} \partial_1 g\left(-N^{-1/\beta}\la,0\right)\pi(V<-\la)d\la
\\
&&
\stackrel{\eqref{eq:tails:1}}{=}\int_{-\infty}^{+\infty} \partial_1 g\left(\la,0\right)\frac{c_{\beta}(\la)}{|\la|^{\beta}}d\la.
\end{eqnarray*}
Likewise
$$
\lim_{N\to+\infty}\tilde{\cal K}_N^{(2)}=\int_{0}^{M} \partial_2 g\left(0,\la\right)\frac{c_{\al}(\la)}{\la^{\al}}d\la
$$
and \eqref{021701b} follows.

 {\em The case \eqref{sing-2} holds.} 
 We shall need the following lemma.
\begin{lemma}
\label{lm020402b}
Suppose that \eqref{sing-2} holds. Then, there exists a constant $C_*$ such that
\begin{equation}
\label{sing-2a}
\pi\left[|\tau-\rho\circ \chi|\ge \lambda\right] \le \frac{C_*}{\lambda^{\gamma}}.
\end{equation}
The exponent $\gamma$ is the same as in \eqref{sing-2}.
\end{lemma}
\proof
The left hand side of \eqref{sing-2a} can be estimated by
\begin{eqnarray}
\label{041801}
&&
\pi\left[|\tau-\rho\circ V|\ge \lambda/2\right]+\pi\left[|\rho\circ V-\rho\circ \chi|\ge \lambda/2\right].
\end{eqnarray}
The first term can be estimated directly from \eqref{sing-2}.
 When $\beta\le \al$ we have
$|\rho(\lambda_1)-\rho(\lambda_2)|\le \rho(\lambda_1-\lambda_2)$ for all
$\lambda_1,\lambda_2\in\bbR$. The second term in \eqref{041801}
can be estimated from \eqref{062202} by
\begin{eqnarray}
\label{041801b}
&&
\pi\left[|P \chi|\ge (\lambda/2)^{\al/\beta}\right]\le \frac{C\|P\chi\|_{L^{2}(\pi)}^{2}}{\lambda^{2\al/\beta}}
\end{eqnarray}
and \eqref{sing-2a} follows for $\beta\in(0,2)$.

When, on the other hand $\beta>\al$ the second term in \eqref{041801}
can be estimated by
\begin{eqnarray}
\label{041802}
&&
\pi\left[\left|\max\{e'\circ V,e'\circ \chi\}\right|\left|(V-\chi)\right|\ge \la\right] \\
&&
\le \pi\left[|e'\circ V P\chi|\ge \lambda/2\right]+\pi\left[|e'\circ
  \chi P\chi|\ge \lambda/2\right].\nonumber
\end{eqnarray}
To estimate the first term on the right hand side we recall that according to Young's
inequality $\lambda_1\lambda_2\le \lambda_1^p/p+ \lambda_2^q/q$ for
any $\lambda_1,\lambda_2>0$ and $p,q>0$ such that $p^{-1}+q^{-1}=1$.
Choose $p$ such that $p_1:=p(\beta/\al-1)<\beta/\al$ and
$(\beta+2)/(2\al)>q>\beta/\al$.
The first term in \eqref{041802} can be estimated by
\begin{eqnarray}
\label{041803}
&&
\pi\left[|V|\ge C_1\lambda^{1/p_1}\right]+\pi\left[|P\chi|\ge C_2\lambda^{1/q}\right]
\end{eqnarray}
for some constants $C_1,C_2>0$ independent of $\lambda$. The first term
can be estimated by $C\lambda^{-\beta/p_1}$, while the second by
$
C\lambda^{-2/q}\|P\chi\|_{L^{2}(\pi)}.
$
These together yield the desired bound on the first term in
\eqref{041802}. The second term can be dealt with similarly.
\qed

Let 
$$
A_N:=[|V|\ge a_*N^{1/\beta}/2,\mbox{ or } \tau^{(N)}\ge a_*N^{1/\al}/2]
$$ 
and for some $\gamma>\kappa>\al$  we let
$
B_N:=[|\rho\circ V-\tau|\ge N^{1/\kappa}].
$ 
Observe that from \eqref{sing-3}, \eqref{sing-2} and Lemma \ref{lm020402b}
\begin{equation}
\label{022005}
\pi(A_N)\le \frac{C}{N}\quad\mbox{ and }\quad\pi(B_N)\le \frac{C}{N^{\gamma/\kappa}}.
\end{equation}

Let also $\rho^{(N)}(x):= \rho(x)1_{[\rho<M_N]}(x)$. Define
$$
\tilde z^{(N)}(x):=(N^{-1/\beta} V(x),N^{-1/\al}
 \rho^{(N)}\circ V(x)).
$$
Note that $ z^{(N)}(x)=\tilde z^{(N)}(x)+r^{(N)}(x)$, where
$$
r^{(N)}(x):=(0,N^{-1/\al} [ \tau^{(N)} (x)- \rho^{(N)}\circ V(x)]).
$$
Note that $z^{(N)}(x)$ lies outside the support of $g$ on
$A_N^c$. Therefore, the expression under the limit in \eqref{021701b} can be written as
\begin{eqnarray*}
&&
N\int g\left(z^{(N)}(x)\right)1_{A_N}\pi(dx)={\cal I}_N+{\cal J}_N,
\end{eqnarray*}
where
\begin{eqnarray*}
&&
{\cal I}_N:=N\int g\left(z^{(N)}(x)\right) 1_{A_N}  1_{B_N}\pi(dx),\quad
{\cal J}_N:=N\int g\left(z^{(N)}(x)\right)  1_{A_N} 1_{B_N^c}\pi(dx).
\end{eqnarray*}
Note that
\begin{eqnarray*}
&&
{\cal I}_N\le N\|g\|_\infty\pi(B_N)\stackrel{\eqref{022005}}{\le}  CN^{1-\gamma/\kappa}\|g\|_\infty\to0,
\end{eqnarray*}
as $N\to+\infty$.
Finally, ${\cal J}_N={\cal J}_N^{(1)}+{\cal J}_N^{(2)}$, where
\begin{eqnarray*}
&&
{\cal J}_N^{(1)}:=N\int g\left(\tilde z^{(N)}(x)\right)1_{A_N}1_{B_N^c}\pi(dx),\\
&&
{\cal J}_N^{(2)}:=N\int\int_0^1 \nabla g\left(\tilde z^{(N)}(x)+\lambda
  r^{(N)}(x)\right) \cdot r^{(N)}(x)1_{A_N}1_{B_N^c}\pi(dx)d\lambda.
\end{eqnarray*}
Given $\delta>0$  we choose $N_0$ so that for $N\ge
N_0$ we have $N^{1/\kappa}<\delta N^{1/\al}$. Let
\begin{eqnarray*}
&&
C_N^{(1)}:=[\rho\circ V<M N^{1/\al},(M+\delta) N^{1/\al}\le \tau],\\
&&
C_N^{(2)}:=[\tau<M N^{1/\al},(M+\delta) N^{1/\al}\le
\rho\circ V]\\
\end{eqnarray*}
 and  
$
D_N:=[N^{-1/\al}\tau\in(M-\delta,M+\delta)].
$
Note that (recall $M_N:=MN^{1/\al}$)
$$
B_N^c\cap(C_N^{(1)}\cup C_N^{(2)}\cup D_N)^c\subset E_N:=(B_N^c\cap[\rho\circ V\le M_N,\tau\le M_N])\cup (B_N^c\cap[\rho\circ V\ge M_N,\tau\ge M_N]).
$$
We have
\begin{eqnarray*}
&&
{\cal J}_N^{(2)}\le N\|\nabla g\|_{\infty}\int |r^{(N)}(x)|1_{A_N}1_{B_N^c}\pi(dx)\\
&&
\le CN^{1+1/\kappa-1/\al}\|\nabla g\|_{\infty}\pi(A_N) +N\|\nabla g\|_{\infty}[\pi(C_N^{(1)})+\pi(C_N^{(2)})+\pi(D_N)].
\end{eqnarray*}
The first term on the right hand side comes from the estimate 
$|r^{(N)}(x)|\le N^{1/\kappa-1/\al}$ that holds on
$E_N$. 
The remaining terms can be estimated by 
$$
c_\al(1+o(1))\|\nabla g\|_{\infty}[(M-\delta)^{-\al}-(M+\delta)^{-\al}],
$$
where $o(1)\to0$, as $N\to+\infty$. This expression can be made arbitrarily small when $\delta>0$ 
is chosen sufficiently small.

Concerning the term ${\cal J}_N^{(1)}$, we can repeat the above argument
and justify in that way that it is equal, up to a term of order $o(1)$, to
\begin{eqnarray*}
&&
\tilde {\cal J}_N^{(1)}:=N\int g\left(\tilde z^{(N)}(x)\right)\pi(dx)\\
&&
=N\int\int_0^{+\infty}\frac{d}{d\lambda} g\left(N^{-1/\beta}\lambda,N^{-1/\al}\rho^{(N)}(\la)\right)1_{[0<\la<V(x)]}\pi(dx)d\lambda
\nonumber\\
&&
-
N\int\int_{-\infty}^{0}\frac{d}{d\lambda} g\left(N^{-1/\beta}\lambda,N^{-1/\al}\rho^{(N)}(\la)\right)1_{[0>\la>V(x)]}\pi(dx)d\lambda.\nonumber
\end{eqnarray*}
Consider the first term on the utmost right hand side of the above expression. Integrating over $x$ we obtain that it equals
\begin{eqnarray*}
&&
N\int_0^{+\infty}\frac{d}{d\lambda} g\left(N^{-1/\beta}\lambda,N^{-1/\al}\rho^{(N)}(\la)\right)\pi[\la<V]d\lambda.
\end{eqnarray*}
Changing variables $\lambda':=\lambda N^{1/\beta}$ and letting $N\to+\infty$ we obtain that the limit equals
\begin{eqnarray*}
&&
c_{\beta}^+\int_0^{+\infty}\frac{d}{d\lambda} g\left(\lambda,1_{[0,M]}(\la)\rho(\la)\right)\lambda^{-\beta}d\lambda=\int_0^{+\infty} g\left(\lambda,1_{[0,M]}(\la)\rho(\la)\right)\nu_{\beta}(d\lambda).
\end{eqnarray*}
The limit for the second term is computed in the same way and we obtain
\eqref{021701b}.
\qed

\subsubsection*{The proof of \eqref{eq:16bc}}
Since $\{X_n,\,n\ge0\}$ is stationary and Markovian it suffices only to show that for $g$ as in the statement of Theorem \ref{prop010805} we have
\begin{equation}
   \label{eq:16zz}
   \lim_{N\to+\infty}N\mathbb E\left\{\mathbb E\left[g\left(
        Z_{1,N}\right)\Big|{\cal G}_0\right]\right\}^2=0.
  \end{equation}
Let  $\delta:=1/2$dist$(0,$supp $g)>0$. We can estimate the expression under the limit by
\begin{eqnarray}
\label{022702}
&&  N\|g\|_{\infty}
\mathbb E\left\{\mathbb P\left[
        |Z_{1,N}|\ge \delta\Big|{\cal G}_0\right]\right\}^2        \le  2\|g\|_{\infty}\left\{N\mathbb E\left\{\mathbb P\left[
        |R_0(X_1,X_0)|\ge \delta N^{1/\beta}/2\Big|{\cal G}_0\right]\right\}^2 
       \right. \nonumber\\
       &&
       \left.+N\mathbb E\left\{\mathbb P\left[
        |\tau^{(N)}(X_1)|\ge \delta N^{1/\al}/2\Big|{\cal G}_0\right]\right\}^2\right\}.
  \end{eqnarray}
The first term on the right hand side of \eqref{022702} can be estimated by
\begin{eqnarray*}
&&  2\|g\|_{\infty}N^{1-2/\beta}\mathbb E\left\{\mathbb E\left[
        |R_0(X_1,X_0)|\Big|{\cal G}_0\right]\right\}^2\le CN^{1-2/\beta}\|P|\chi|\|_{L^2(\pi)}^2
         \end{eqnarray*}
         for some constant $C>0$ independent of $N$. The expression on the right hand side tends to $0$
as $N\to+\infty$, thanks to \eqref{062202} and $\beta\in(1,2)$.

The second term on the right hand side of \eqref{022702} can be estimated, using \eqref{022602}, by
\begin{eqnarray*}
&&  2\|g\|_{\infty}N^{1-2/\al}\mathbb E\left\{\mathbb E\left[
        \tau^{(N)}(X_1)\Big|{\cal G}_0\right]\right\}^2\\
        &&
        \le 2 \|g\|_{\infty}N^{1-2/\al}\|P\tau^{N)}\|_{L^2(\pi)}^2\le
        CN^{1-2/\al}N^{2(1/\al-1)}=\frac{C}{N}
         \end{eqnarray*}
         for some constants $C>0$ independent of $N$.
         
\subsubsection*{The proof of \eqref{eq:16bcde}}         
  Using again stationarity and Doob's inequality for the martingale  $\{Z_t^{(\Delta,N)},\,t\ge0\}$ we obtain that    the expression under the limit in \eqref{eq:16bcde} can be estimated from above by
   $$
   CTN\sum_{j=1}^2\left\{\bbE \left[Z_{1,N}^{(j)}1_{[|Z_{1,N}^{(i)}|\le \Delta_i,i=1,2]}\right]^2-\bbE\left[\bbE\left[Z_{1,N}^{(j)}1_{[|Z_{1,N}^{(i)}|\le \Delta_i,i=1,2]}\Big|{\cal G}_{0,N}\right]\right]^2\right\}.
   $$  
   The term corresponding to $j=1$ can be estimated by
   \begin{eqnarray}
   \label{052702}
&&   C_1TN^{1-2/\beta}\int \chi^2(x)1_{[|\chi(x)|\le \Delta_1 N^{1/\beta}]}\pi(dx)
\\
&&
\le C_2TN^{1-2/\beta}\int_0^{\Delta_1 N^{1/\beta}}\la \pi(|\chi(x)|\ge \la)d\la
= C_2TN\int_0^{\Delta_1 }\la \pi(|\chi(x)|\ge N^{1/\beta}\la)d\la\nonumber
   \end{eqnarray}
   for some constants $C_1,C_2>0$.
   Using the fact that $ \pi(|\chi(x)|\ge N^{1/\beta}\la)\le C_3N^{-1}\la^{-\beta}$ for some constant $C_3>0$ we conclude that the utmost right hand side of \eqref{052702} can be estimated by
   $$
   CT\int_0^{\Delta_1 }\la^{1-\beta} d\la=C'T\Delta_1^{2-\beta }
   $$
   for some constants $C,C'>0$. A similar estimate can be obtained also for $j=2$. Estimate \eqref{eq:16bcde} then follows.

\section{The proof of Theorem \ref{joint-conv}}

\label{appB}

 In this section  we retain the notation from the previous section. Assume furthermore that $\beta=2$. It is well known that under the hypotheses made in Section \ref{markov-chain}, the components of processes $\{X_t^{(N)}:=(S_t^{(N)},T_t^{(N)}),\,t\ge0\}$  converge in ${\cal D}$, see e.g. \cite{jakola,rosenblatt}, to a Brownian motion and an $\alpha$-stable subordinator process. This in turn implies
tightness of the laws of $\{X_t^{(N)},\,t\ge0\}$  over ${\cal D}\times{\cal D}$ equipped with the uniform and $J_1$ topologies, respectively. To finish the proof of the theorem in question we only need to show the weak convergence of finite dimensional distributions of $\{X_t^{(N)},\,t\ge0\}$
to the respective finite dimensional distribution of  a L\'evy process
$\{( S_t, T_t),\,t\ge0\}$ with the exponent 
\begin{eqnarray}
\label{020103}
\psi(\xi_1,\xi_2):=\si^2\xi_1^2/2+
\int_{0}^{+\infty}(e^{i\xi_2 \la}
-1)\nu_\al(d\la)
\end{eqnarray}
for some $\si\ge0$.
By a well known  Cram\'er-Wold device, see Theorem 9.5, p. 147 of \cite{durrett}, it suffices to consider only one dimensional distributions
i.e. to prove that for $(\xi_1,\xi_2)\in\bbR^2$ and $t\ge0$
\begin{equation}
\label{112702}
\xi_1  S^{(N)}_{t}+\xi_2 T^{(N)}_{t}\Rightarrow\xi_1  S_{t}+\xi_2 T_t,\quad\mbox{as }N\to+\infty
\end{equation}

 For a fixed $M>0$ define $\bar Z^{(2)}_{N,t}(M):=\sum_{n=0}^{[Nt]-1}\bbE[\tau^{(N)}(X_n)|{\cal G}_{n-1,N}]$ with $\tau^{(N)}(x)$ given by \eqref{010103}. Using \eqref{eq:16bc} one can argue, as in the proof of (4.3) of \cite{durrett}, that
 \begin{equation}
  \label{102702}
  \lim_{N\to+\infty}\bar Z^{(2)}_{N,t}(M)=
    t\int_0^M \la\nu_{\al}(d\la)
\end{equation}
in probability. Since \eqref{012702} is also in force for $\beta=2$, thanks to \eqref{102702},
 it suffices only to show that for a fixed $M>0$
 \begin{equation}
\label{112702b}
\xi_1  Z^{(1)}_{N,t}+\xi_2 \tilde Z^{(2)}_{N,t}(M)\Rightarrow\xi_1  S_{t}+\xi_2 \tilde T_t(M),\quad\mbox{as }N\to+\infty,
\end{equation}
where  $\tilde Z^{(2)}_{N,t}(M)=Z^{(2)}_{N,t}(M)-\bar Z^{(2)}_{N,t}(M)$, $Z^{(1)}_{N,t}$, $Z^{(2)}_{N,t}(M)$ are given by
 \eqref{122702}
and $\{(S_t,\tilde T_t(M)),\,t\ge0\}$ is a L\'evy process with the exponent.
 \begin{eqnarray}
\label{020103b}
\psi(\xi_1,\xi_2):=\si^2\xi_1^2/2+
\int_{0}^M(e^{i\xi_2 \la}
-1-i\xi_2\la)\nu_\al(d\la)
\end{eqnarray}
 for some $\si>0$. In what follows we omit writing $M$ in the notation of processes.
Note that  $\{(Z^{(1)}_{N,t}, Z^{(2)}_{N,t}),\,t\ge0\}$ is a martingale with the increments given by
$$
\tilde Z_{n,N}=(\tilde Z^{(1)}_{n,N},\tilde Z^{(2)}_{n,N}):=\left(\frac{1}{N^{1/2}}R_0(X_n,X_{n-1}),\frac{1}{N^{1/\al}}\left\{\tau^{(N)}(X_n)-\bbE[\tau^{(N)}(X_n)|{\cal G}_{n-1,N}]\right\}\right) . 
$$
We can use therefore the results of  \cite{brown-eagleson}. According to Theorem 1 of ibid. in order to show \eqref{112702b} it suffices to prove the following result.
\begin{prop}
\label{lim-law}
Suppose that $M_{-1,N}:=0$ and
$$
M_{n,N}:=\frac{\xi_1}{N^{1/2}}R_0(X_n,X_{n-1})+\frac{\xi_2}{N^{1/\al}}\left\{\tau^{(N)}(X_n)-\bbE[\tau^{(N)}(X_n)|{\cal G}_{n-1,N}]\right\} 
$$
for $n\ge0$.
Then, for any $a<b$ we have
\begin{equation}
\label{fundamental}
\lim_{N\to+\infty}\sum_{n=1}^{[Nt]-1}\mathbb E\left[M_{n,N}^21_{[a<M_{n,N}<b]}\left|\right.{\cal G}_{n-1,N}\right]=t[G_{\xi_1,\xi_2}(b)-G_{\xi_1,\xi_2}(a)]
\end{equation}
in probability, where
the function $G_{\theta_1,\theta_2}(\cdot)$ is given by $G_{\xi_1,\xi_2}(\la)=0$, $\la<0$ and
$$
G_{\xi_1,\xi_2}(\la)=\si^2\xi_1^2+c_\al |\xi_2|^\al\la^{1-\al}.
$$
for $\la>0$ and some $\si>0$.
\end{prop}
\proof 
The proof is done in two steps. First we
prove that for any interval $(a,b)$ that does not contain $0$ and any
$C^\infty$ function $g$ supported in that interval we have
\begin{equation}
\label{case-1}
\lim_{N\to+\infty}\sum_{n=0}^{[Nt]-1}\mathbb E\left[g(M_{n,N})\left|\right.{\cal G}_{n-1,N}\right]=t|\xi_2|^{\al}\int_0^{M}g(\lambda)\nu_{\al}(d\la)
\end{equation}
in probability. The argument  is essentially a repetition of the proof of \eqref{eq:16bc} done in the previous section so we omit it.
Next, we show that for any $c>0$
\begin{equation}
\label{case-2}
\limsup_{N\to+\infty}\bbE\left|\sum_{n=0}^{[Nt]-1}\mathbb
  E\left[M_{n,N}^21_{[|M_{n,N}|<c]}\left|\right.{\cal
      G}_{n-1}\right]-\frac12\si^2t\xi_1^{2}\right|= h(c),
\end{equation}
where $\lim_{c\to0+}h(c)=0$.
For an arbitrary interval $(a,b)$, where $a<0<b$ we divide it into
a  sum of three disjoint intervals $(a,-c)$, $(-c,c)$ and $(c,b)$,
where $0<c<\min[-a,b]$ and conclude using the above results and a
standard approximation argument that 
\begin{equation*}
\label{case-3}
\limsup_{N\to+\infty}\bbE\left|\sum_{n=0}^{[Nt]-1}\mathbb
  E\left[M_{n,N}^21_{[a<M_{n,N}<b]}\left|\right.{\cal
      G}_{n-1,N}\right]-t[G_{\xi_1,\xi_2}(b)-G_{\xi_1,\xi_2}(a)]\right|=0.
\end{equation*}
\subsection*{The proof of  (\ref{case-2})} 
Suppose that $c>0$ is arbitrary. We consider only the case when both $\xi_1,\xi_2\not=0$. The other cases can be done  adjusting (and simplifying) the argument.
Note that
\begin{eqnarray*}
&&\sum_{n=0}^{[Nt]-1}\mathbb E\left[(M_{n,N})^21_{[|M_{n,N}|<c]}\left|\right.{\cal G}_{n-1,N}\right]
=\xi^2_1\sum_{n=0}^{[Nt]-1}\mathbb E\left[(\tilde Z_{n,N}^{(1)})^21_{[|M_{n,N}|<c]}\left|\right.{\cal G}_{n-1,N}\right]
\\
&&
+\xi^2_2\sum_{n=0}^{[Nt]-1}\mathbb E\left[(\tilde Z_{n,N}^{(2)})^21_{[|M_{n,N}|<c]}\left|\right.{\cal G}_{n-1,N}\right]
+2\xi_1\xi_2
\sum_{n=0}^{[Nt]-1}\mathbb E\left[\tilde Z_{n,N}^{(1)}\tilde Z_{n,N}^{(2)}1_{[|M_{n,N}|<c]}\left|\right.{\cal G}_{n-1,N
}\right].
\end{eqnarray*}
Denote the terms appearing on the right hand side by $U_N$, $V_N$ and $W_N$.
For  an appropriate constant $C>0$ we have
\begin{eqnarray*}
&&
\bbE|W_N|\le C
\left\{N\left\{\mathbb E[[\xi_1 \tilde Z_{1,N}^{(1)}]^2
]\right\}^{1/2}\left\{\mathbb E[[\xi_2\tilde  Z_{1,N}^{(2)}]^2
,\,|\xi_2 \tilde Z_{1,N}^{(2)}|<10c\right\}^{1/2}\right.\\
&&
+\left.N|\xi_1\xi_2|\mathbb E\left[|\tilde Z_{1,N}^{(1)}\tilde Z_{1,N}^{(2)}|,\,|\xi_1 \tilde Z_{1,N}^{(1)}|>9c, |\xi_2 \tilde Z_{1,N}^{(2)}|\ge 10c
\right]\right\}
\end{eqnarray*}
Denote the  first and  second term appearing in the braces on the right hand side  by $W_N^{(1)}$ and $W_N^{(2)}$ respectively. 
We have
\begin{eqnarray*}
\mathbb P[|\xi_2 \tilde Z_{n,N}^{(2)}|> \la]
\le \pi[\tau > N^{1/\al}\la/(2|\xi_2| )] + \pi[
P \tau^{(N)}> N^{1/\al}\la/(2 |\xi_2|)].\nonumber
\end{eqnarray*}
The first term on the right hand side is clearly less than, or equal to $ CN^{-1}\la^{-\al}$
for  all $\la>0$, $N\ge1$ and a certain constant $C>0$, independent of
$n\ge0$.
By the Markov inequality and \eqref{022602} the second term can be estimated by
$
\la^{-1}N^{-1/\al}\|P\tau\|_{L^2(\pi)}\le
C(\la N)^{-1}.
$ 
One can deduce that
\begin{equation}
\label{031702}
\mathbb P[|\xi_2 \tilde Z_{n,N}^{(2)}|> \la]\le C(N\la)^{-1}(1+\la^{1-\al}) 
\end{equation}
for  all $\la>0$, $N\ge1$ and a certain constant $C>0$, independent of
$n,N\ge0$.
Using \eqref{031702} and an elementary estimate 
$
\left\{\mathbb E[[\xi_1 \tilde Z_{1,N}^{(1)}]^2
]\right\}^{1/2}\le CN^{-1/2}\|\chi\|_{L^2(\pi)}
$
 we obtain
\begin{equation*}
\label{021802}
W_N^{(1)}\le C\|\chi\|_{L^2(\pi)}\left[\int_0^{10c}(1+\la^{1-\al})d\la \right]^{1/2}\le C'\|\chi\|_{L^2(\pi)}[c(1+c^{1-\al})]^{1/2}
\end{equation*}
for some constants $C,C'>0$.

On the other hand, 
 using Chebyshev's inequality we get
\begin{eqnarray}
\label{041802b}
&&
\mathbb P[|\xi_1\tilde Z_{n,N}^{(1)}|> \la]\le \frac{C\|\chi\|_{L^2(\pi)}^2}{N\la^2}
\end{eqnarray}
for all $\la>0$. The constant $C>0$ appearing here does not depend on $N,n$ and $\lambda$. Thus,
for some constants $C,C'>0$, we have
\begin{eqnarray*}
&&W_N^{(2)}\le C NM\mathbb E\left[|\xi_1\tilde Z_{1,N}^{(1)}|,\,|\xi_1 \tilde Z_{1,N}^{(1)}|>9c
\right]\\
&&
\le C NM\left\{
\mathbb E\left[[\xi_1\tilde Z_{1,N}^{(1)}]^2,\,|\xi_1\tilde Z_{1,N}^{(1)}|>9c
\right]\right\}^{1/2} 
\mathbb P^{1/2}\left[|\xi_1\tilde Z_{1,N}^{(1)}|>9c
\right]\nonumber\\
&&
\stackrel{\eqref{041802b}}{\le} C'\left\{
\mathbb E\left[R^2_0(X_1,X_0),\,|\xi_1 \tilde Z_{1,N}^{(1)}|>9c
\right]\right\}^{1/2} \to0,\nonumber
\end{eqnarray*}
 as $N\to+\infty$. 
We have proved therefore that
$
\limsup_{N\to+\infty}\bbE|W_N|\le C [c(1+c^{1-\al})]^{1/2},
$
where $c>0$ can be chosen arbitrarily. Thus, $\lim_{N\to+\infty}\bbE|W_N|=0$.

Note that
\begin{eqnarray}
\label{061802b}
&&
\mathbb E|V_N|\le CN\mathbb E
[[\xi_2 \tilde Z_{n,N}^{(2)}]^2,\,|\xi_2\tilde Z_{n,N}^{(2)}|<10c]+
\\
&&
+CN^{1-2/\al}\mathbb E[
[ \tau^{(N)}(X_0)]^2,\,|\xi_1\tilde Z_{n,N}^{(1)}|>9c,\,|\xi_2\tilde Z_{n,N}^{(2)}|\ge 10c]\nonumber
\end{eqnarray}
for some constant $C>0$. Denote the first and the second terms on the right hand side of \eqref{061802b} by $V_N^{(1)}$, $V_N^{(2)}$ respectively.
Estimating in the same way as in  \eqref{031702} we deduce
$
V_N^{(1)}\le  C [c(1+c^{1-\al})]^{1/2}
$
for some constant $C>0$. This term can be made arbitrarily small by choosing a sufficiently small $c>0$.
On the other hand, from Chebyshev's inequality for some constants $C,C'>0$ we have
\begin{eqnarray*}
&&
V_N^{(2)}\le CNM^2 
\mathbb P[|\xi_1 R_0(X_1,X_0)|\ge 9cN^{1/2}]\\
&&
\le  C'M^2 \mathbb E[[ \xi_1R_0(X_1,X_0)]^2, |\xi_1R_0(X_1,X_0)|\ge 9cN^{1/2}]\to0,
\end{eqnarray*}
both a.s. and in the $L^1$ sense,
as $N\to+\infty$.
 Finally,
we can write that
$
U_N=\widehat U_N-\bar U_N,
$
where
$$
\widehat U_N:=\xi^2_1\sum_{n=0}^{[Nt]-1}\mathbb E\left[(\tilde Z_{n,N}^{(1)})^2\left|\right.{\cal G}_{n-1,N}\right]
\quad\mbox{and}\quad
\bar U_N:=\xi^2_1\sum_{n=1}^{[Nt]}\mathbb E\left[(\tilde Z_{n,N}^{(1)})^21_{[|\tilde Z_{n,N}|>c]}\left|\right.{\cal G}_{n-1,N}\right].
$$
By the ergodic theorem
$$
\widehat U_N=\frac{\xi^2_1}{N}\sum_{n=0}^{[Nt]-1}\left[PV^2(X_{n-1}) +P\chi^2(X_{n-1})-(P\chi)^2(X_{n-1})\right]\to \frac12\sigma^2\xi^2_1 t,\quad\mbox{as }N\to+\infty,
$$
both a.s. and in the $L^1$ sense. Here 
$
\sigma^2:=2\left(\|V\|_{L^2(\pi)}^2+\|\chi\|_{L^2(\pi)}^2-\|P\chi\|_{L^2(\pi)}^2\right)
.$
Using stationarity we deduce that for a certain constant $C>0$
\begin{eqnarray*}
\mathbb E|\bar U_N|\le 
C\xi^2_1\mathbb E\left[R_0^2(X_1,X_0),|\xi_1R_0(X_1,X_0)|>cN^{1/2}/2,\mbox{ or }
|\xi_2\tilde Z_{1,N}^{(2)}|>c/2\right]\to 0,
\end{eqnarray*}
as $N\to+\infty$.
The convergence follows from the $L^2$-integrability of $R_0(X_1,X_0)$ and \eqref{031702}. 
\qed

\end{appendix}

\end{document}